\documentclass{book}

\usepackage{mathtools,amsfonts,amssymb,amsthm}
\usepackage{color,graphicx,cite,float,url}
\usepackage[center]{caption}

\title{Notes of the Summer School \\ on Finite Set Statistics}
\author{}
\date{}

\newtheorem{definition}{Definition}
\newtheorem{example}{Example}
\newtheorem{remark}{Remark}

\newtheorem{myexercise}{Exercise}
\newtheorem{mytheorem}{Theorem}
\newtheorem{mylemma}{Lemma}
\newtheorem{myproof}{Proof}
\newtheorem{myrule}{Rule}

\def\eqn#1{\begin{equation}#1\end{equation}}

\def\xline#1{\begin{multline}#1\end{multline}}
\def\eqnalign#1{\begin{align}#1\end{align}}
\def\subeqn#1#2{\begin{subequations}\label{#1}\begin{align}#2\end{align}\end{subequations}}
\def\itemName#1{\item $[${\em#1}$\,]$}
\def\ind#1{\mathbf{1}_{#1}}
\def\dd{\mathrm{d}}
\def\vphi{\varphi}
\def\ii{\mathrm{i}}
\def\lrarrow{\leftrightarrow}
\def\target{\text{target}}
\def\birth{\text{birth}}
\def\clutter{\text{clutter}}
\def\obs{\text{obs}}

\def\BorelMeasSpace#1{(#1,\mathcal{B}(#1))}

\def\R{\mathbb{R}}
\def\var{\mathop{\rm var}} 

\newenvironment{hint}{~\\ \textit{(Hint)}\itshape}{~\\ \\}

\newcommand{\Bcal}{\mathcal{B}} 

\newcommand{\Xcal}{\mathcal{X}} 
\newcommand{\Nset}{\mathbb{N}} 
\newcommand{\ds}{\displaystyle} 
\newcommand{\Exp}{\mathbb{E}} 

\begin{document}

\maketitle

\chapter[Measure and probability theory]{Measure and probability theory for multi-object estimation}
\label{chapt:MeasureAndProba}

The main objective in this lecture is to give an overview of the concepts and tools of probability theory and measure theory that are useful for multi-object estimation. Probability and measure theory are usually described in a very formal way, and the underlying ideas behind the introduced concepts might be difficult to grasp.

In spite of the limited duration of this lecture, both measure and probability theory will be discussed. Indeed,
\begin{quote}
``Probability theory and measure theory have become so intertwined that they seem to many mathematicians of our generation to be two aspects of the same subject''\footnote{M. Adams and V. Guillemin. Measure Theory and Probability, Wadsworth \& Brooks, Monterey, California, 1986.}.
\end{quote}

This lecture is organised as follows. Starting from simple consideration about a random experiment, questions will be raised and answered successively until the concept of point process is reached and described. Indeed, point processes are useful to describe random experiments that are sophisticated enough to deal with multi-object estimation. Some details will be omitted when they would provide only little insight about the general concept.

Consider a space $\mathcal{X}$ (nice enough, e.g.\ $\mathbb{R}^d$ for any $d > 0$) where we want to describe a random experiment. The space $\mathcal{X}$ is called the outcome space since it contains all the possible outcomes regarding the random experiment of interest, as illustrated in the following example. 

\begin{example}[Where is the fly?]
\label{ex:flyFirst}
Let's try a thought experiment. Consider that there is a fly in a room, and that its motion is random. Let $S \subset \mathbb{R}^3$ be the volume of the room. An outcome for the position of the fly is any point $x \in S$ so that the space of all outcomes is $S$. Now, if you stop observing the fly and if, after a couple of seconds, you try to predict its position, you are likely to say ``it should be around there'' rather than ``it is exactly there''. This simple thought experiment shows that outcomes are not so useful to describe random experiments; the introduction of another concept is required.
\end{example}

\subsection*{Q1: How to put mathematical sense on ``around there''?}

\begin{figure}[h]
\centering
\includegraphics[width=300px]{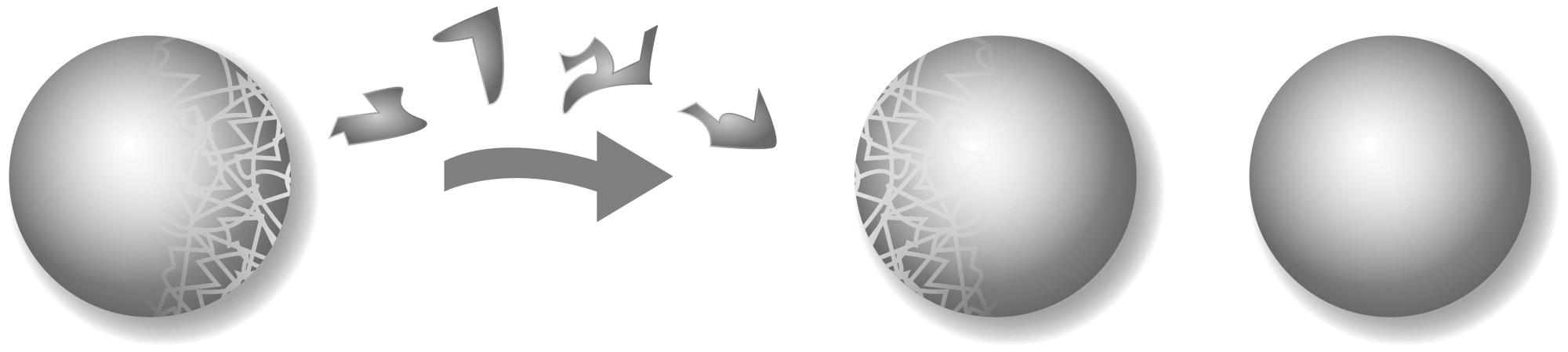}
\caption{The Banach-Tarski paradox}
\label{fig:BanachTarski}
\end{figure}

Clearly, the concept of subset is required. However, it is not desirable to consider any subset of the space $\mathcal{X}$, i.e.\ any element in the power set of $\mathcal{X}$. Indeed, some subsets have very peculiar properties that are far from being intuitive. The Banach-Tarski paradox illustrated in Figure \ref{fig:BanachTarski} shows that, when considering any subset of a ball, one can build two balls out of one, all of them being identical.

It is therefore required to only consider the subsets of $\mathcal{X}$ that have nice properties, and we call them the measurable subsets of $\mathcal{X}$. We additionally want to consider specific classes of measurable sets, with the following properties
\begin{enumerate}
\item closed under complementation and countable union,
\item contains the empty set $\emptyset$.
\end{enumerate}
Any class $\mathcal{S}$ of subsets of $\mathcal{X}$ with these properties is called a $\sigma$-algebra on $\mathcal{X}$. The space $(\mathcal{X},\mathcal{S})$ is called a measurable space. 

\subsubsection*{Q1': Do we want any $\sigma$-algebra?}

No, because some $\sigma$-algebras are too simple to be useful, for instance the minimal/trivial $\sigma$-algebra is the set $\{\emptyset,\mathcal{X}\}$. A useful $\sigma$-algebra is called the Borel $\sigma$-algebra and is defined as the smallest $\sigma$-algebra containing all the open sets in $\mathcal{X}$, denoted $\mathcal{B}(\mathcal{X})$.

\subsection*{Q2: How can we use/characterise these measurable subsets?}

To use or characterise measurable subsets, we need a function taking measurable subsets as argument and associating a real number (possibly infinite).

\begin{definition}[Measure] Let $(\mathcal{X},\mathcal{S})$ be a measurable space, a function
\eqn{
\mu:\mathcal{S} \to \mathbb{R}\cup\{+\infty\}
}
is called a measure if:
\begin{enumerate}
\item \label{it:measPositive} $\forall B \in \mathcal{S}$: $\mu(B) \geq 0$,
\item \label{it:measZeroEmptySet} $\mu(\emptyset) = 0$,
\item \label{it:measSigmaAdd} $\forall B_1, B_2, \ldots \in \mathcal{S}$, $B_i \cap B_j = \emptyset$, $\mu(\cup B_i) = \sum \mu(B_i)$,
\end{enumerate}
Additionally, $\mu$ is called a probability measure if 
\eqn{
\label{it:measProba}
\mu(\mathcal{X}) = 1.
}

The space $(\mathcal{X},\mathcal{S},\mu)$ is called a measure space or a probability space if $\mu$ is a probability measure.
\end{definition}

Each of the hypotheses above can be interpreted in the context of probability theory:
\begin{itemize}
\item Hypothesis \ref{it:measPositive} ensures that there is no event with negative probability, which would not be meaningful.
\item Hypothesis \ref{it:measZeroEmptySet} makes sure that the probability of the empty set $\emptyset$ representing the event ``no event'' is null. For instance, the probability for a fair coin to give neither head or tail has to be null.
\item Hypothesis \ref{it:measSigmaAdd} corresponds to the principle that the probability for a fair die to give $2$ or $5$ is the sum of the probability of this two events. Note that, for any measure $\mu$ and any event $B \in \mathcal{S}$,
\eqn{
\mu(B) = \mu(B \cup \emptyset) = \mu(B) + \mu(\emptyset),
}
so that $\mu(\emptyset) = 0$.
\item By convention, the probability of a ``sure event'' is $1$, or ``100\%''. The condition \eqref{it:measProba} ensures that the space $\mathcal{X}$ contains all the possible outcomes so that $\mathcal{X}$ is a sure event.
\end{itemize} 

\begin{remark}[Cartesian product] Let $(\mathcal{X}_1,\mathcal{B}(\mathcal{X}_1),\mu_1)$ and $(\mathcal{X}_2,\mathcal{B}(\mathcal{X}_2),\mu_2)$ be two measure spaces. The Borel $\sigma$-algebra on the product space $\mathcal{X} = \mathcal{X}_1 \times \mathcal{X}_2$ can be expressed as $\mathcal{B}(\mathcal{X}) = \mathcal{B}(\mathcal{X}_1)\otimes\mathcal{B}(\mathcal{X}_2)$ which is a product $\sigma$-algebra.

For any $B_1 \in \mathcal{B}(\mathcal{X}_1)$ and any $B_2 \in \mathcal{B}(\mathcal{X}_2)$, the product measure $\mu = \mu_1 \times \mu_2$ on $(\mathcal{X},\mathcal{B}(\mathcal{X}))$ is such that
\eqn{
\mu(B_1 \times B_2) = \mu_1(B_1) \mu_2(B_2).
}

Measures on product spaces will be useful when we will be dealing with multi-variate random experiments.
\end{remark}

One of the most fundamental and useful example of measures is the Lebesgue measure $\lambda$ on $\mathcal{X} = \mathbb{R}^d$, $d > 0$, endowed with the Borel $\sigma$-algebra $\mathcal{B}(\mathcal{X})$, that gives the ``volume'' $\lambda(B)$ of any measurable subset $B \in \mathcal{B}(\mathcal{X})$.

\subsection*{Q3: Can we only measure subsets?}

It would be nice to be able to measure functions. Let $f$ be a function such that,
\eqn{
\label{eq:simpleFunction}
f = \sum_{i=1}^N a_i \ind{A_i},
}
where $a_i \in \mathbb{R}$, $1 \leq i \leq n$, and where $\ind{A_i}$ is the indicator function of the set $A_i$, i.e.\ $\ind{A_i}(x) = 1$ if $x \in A_i$ and $0$ otherwise. The Lebesgue measure $\lambda(f)$ of the function $f$ can then be defined as
\eqn{
\lambda(f) = \sum_{i=1}^N a_i \lambda(A_i).
}

A very large class of function can be generated by simple functions, and for any of the function $f$ in this class, we can write:
\eqn{
\lambda(f) = \int f \dd \lambda = \int f(x) \lambda(\dd x),
}
which is more general than the Riemann integral. Note that the shorthand notation $\dd x$ is often used to refer to the Lebesgue measure $\lambda(\dd x)$.

The same principle can actually be applied for any measure $\mu$, so that we can write
\eqn{
\mu(f) = \int f \dd \mu = \int f(x) \mu(\dd x),
}
and call it the Lebesgue integral of $f$ with respect to the measure $\mu$. Note that, for any $B \in \mathcal{B}(\mathcal{X})$, $\mu(\ind{B}) = \int \ind{B}(x) \mu(\dd x) = \int_B \mu(\dd x) = \mu(B)$.

Lebesgue integration is not only useful because it is defined for a larger class of functions than the Riemann integral, but also because integration can be carried out on any space as long as it is equipped with a $\sigma$-algebra and a measure. We now have enough tools and concepts to define a random experiment properly. Simple random experiments are dealt with first.

\subsection*{Q4: How to characterise a ``single-variate'' random experiment?}

Here, the focus is on the characterisation of a single individual or entity, e.g.\ one fly, one coin, one die, etc., which justify the name ``single-variate'' random experiment. Since we know how to define a probability measure, the most natural way of characterising a random experiment is to find the right space, say $\mathcal{X}$, to endow it with an adequate $\sigma$-algebra so that a probability measure can be defined.

\begin{example} Consider again the fly from Example \ref{ex:flyFirst}. As stated earlier, the set $S \subset \mathbb{R}^3$ is appropriate to describe the (random) position of the fly. We equip $S$ with the Borel $\sigma$-algebra $\mathcal{B}(S)$ since we are interested in describing the probability for the fly to be in a given area. We can now define a probability measure $P$ so that we have a probability space $(S,\mathcal{B}(S),P)$ to work with.
\end{example}

\begin{figure}[h]
\centering
\def\svgwidth{0.5\columnwidth}
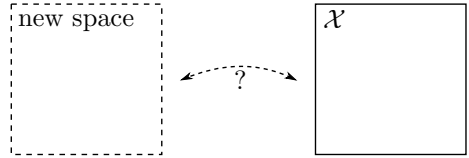
\caption{Are these spaces related?}
\label{fig:outOfTheBlueSpace}
\end{figure}

At this point, we are able to characterise any single variate random experiment. However, it is of interest to be able to relate the events in $\mathcal{X}$ with events in other spaces, as pictured in Figure \ref{fig:outOfTheBlueSpace}. This operation does not seem straightforward because even though the randomness is already well characterised, its source has not been clearly identified. 

\subsubsection*{Q4': Is there a more convenient way of defining a single-variate random experiment?}

\begin{figure}[h]
\centering
\def\svgwidth{0.5\columnwidth}
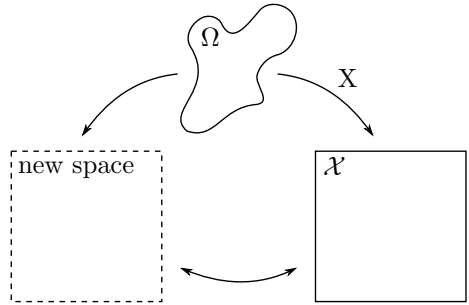
\caption{Proper way of relating spaces.}
\label{fig:WellRelatedSpace}
\end{figure}

To better represent the randomness in one or several spaces, it is necessary to isolate it in a single, possibly complicated probability space $(\Omega,\mathcal{F},\mathbb{P})$ containing everything about the random experiment, see Figure \ref{fig:WellRelatedSpace}.

A random variable is then a nice function $X$ from $(\Omega,\mathcal{F},\mathbb{P})$ to $(\mathcal{X},\mathcal{B}(\mathcal{X}))$. We will come back to what a nice function is in this case after the following example.

\begin{example}[The fly again] If we keep representing the position of our fly in $S$, but at some point we get information about, e.g.\ the age of the fly, then we have another random variable that is characterising the same fly. Clearly, the position and the age of the fly are actually projections of a more complex random experiment where everything about the fly would be characterised. If we represent this latter random experiment by the probability space $(\Omega,\mathcal{F},\mathbb{P})$, then the two random variables, say $X_{\text{position}}$ and $X_{\text{age}}$ are the mappings
\eqn{
X_{\text{position}} : (\Omega,\mathcal{F},\mathbb{P}) \to (S,\mathcal{B}(S)),
}
and
\eqn{
X_{\text{age}} : (\Omega,\mathcal{F},\mathbb{P}) \to (\mathbb{R}^+,\mathcal{B}(\mathbb{R}^+)).
}
\end{example}

\begin{figure}[h]
\centering
\def\svgwidth{0.35\columnwidth}
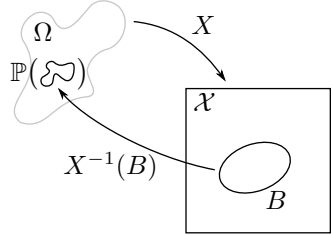
\caption{What is the probability measure of $B$?}
\label{fig:MeasurableFunction}
\end{figure}

So, what are the properties of a function $X$ required to turn it into a proper random variable? Because we are still interested in measuring the probability of events in $\mathcal{X}$, we require that $X$ maps all the events in $\mathcal{X}$ back to events in $\Omega$ so that they can be measured by $\mathbb{P}$, see Figure \ref{fig:MeasurableFunction}. More formally, we require that
\eqn{
\label{eq:measFunction}
X^{-1}(B) = \{ \omega \in \Omega | X(\omega) \in B \} \in \mathcal{F},\quad \forall B \in \mathcal{B}(\mathcal{X}).
}

The functions that satisfy \eqref{eq:measFunction} are called measurable functions. An interesting result is that any non-negative measurable function is the pointwise limit of a sequence of simple functions, as introduced in \eqref{eq:simpleFunction}. The consequence is that any random variable can be integrated in the context of Lebesgue integration. We are now able to define the expectation. However, before doing so, it is useful to check that we can still define a probability measure directly on $\BorelMeasSpace{\mathcal{X}}$.

\subsubsection*{Q4'': What is the probability measure on $(\mathcal{X},\mathcal{B}(\mathcal{X}))$?}

The probability measure on $(\mathcal{X},\mathcal{B}(\mathcal{X}))$ can be defined from $\mathbb{P}$ through the following operation:
\eqn{
P_X(B) = \mathrm{Pr}(X \in B) = \mathbb{P}(X^{-1}(B)), \quad B \in \mathcal{B}(\mathcal{X}).
}
The probability measure $P_X$ is called a pushforward probability measure since it ``pushes'' the probability measure $\mathbb{P}$ ``forward'' into $\mathcal{B}(\mathcal{X})$.

\subsection*{Q5: How to define the expectation?}

Since random variables are measurable functions, we can integrate them with respect to any measure defined on their domain. In particular we would like to sum over all the possible values $X(\omega)$ of $X:(\Omega,\mathcal{F},\mathbb{P}) \to \BorelMeasSpace{\mathcal{X}}$ together with the probability $\mathbb{P}(\dd \omega)$ for the state to be arbitrarily close to $\omega$. Formally, such a quantity can be written as the Lebesgue integral:
\subeqn{eq:expectRV}{
\mathbb{E}[X(\cdot)] & = \int_{\Omega} X \dd\mathbb{P}\\
& = \int_{\Omega} X(\omega) \mathbb{P}(\dd \omega)\\
\label{eq:expectMeas}
& = \int_{\mathcal{X}} x P_X(\dd x),
}
and is called the expectation of $X$, but also the expected value, the mean or the first moment of $X$.

It is worth underlying that the probability space $(\Omega,\mathcal{F},\mathbb{P})$ does not appear explicitly in \eqref{eq:expectMeas}. This space is only used to properly define the random variable $X$, but does not increase the complexity of the expression of the usual concepts like the expectation. Yet, \eqref{eq:expectMeas} is not the usual expression of the expectation since it is based on a probability measure rather than a probability density.

\subsubsection*{Q5': How to relate probability measures and probability densities?}

The first thing to do is to consider a ``small'' measurable subset ``$B = \dd x$'' since we are trying to describe the event at a given point $x$ of the space. However, the mass of any (non-atomic) probability measure in the subset $B$ tends toward $0$ when the size of $B$ tends toward $0$. The quantity of interest is then the ratio $p_X(x)$ between the probability $P_X(\dd x)$ and the size $\lambda(\dd x)$ of $\dd x$, when the latter is arbitrarily small. Informally, we can write
\eqn{
``B = \dd x,\quad p_X(x) = \dfrac{P_X(\dd x)}{\lambda(\dd x)}'',
}
which is well defined when
\eqn{
P_X(B) = 0 \implies \lambda(B) = 0.
}
This condition is formally referred to as the absolute continuity of $P_X$ with respect to $\lambda$ and is denoted $P_X \ll \lambda$.

The function $p_X$ is called the probability density of the random variable $X$ and is formally defined as the function satisfying
\eqn{
P_X(B) = \int_B p_X(x) \lambda(\dd x) = \int_B p_X(x) \dd x.
}

The expectation can now be expressed as
\eqn{
\mathbb{E}[X] = \int x P_X(\dd x) = \int x \, p_X(x) \, \dd x,
}
which is the usual definition of the expectation.

At this point, the single-variate random experiments have been defined and characterised. Many results about them could be written, but our interest is in ``multi-variate'' random experiments which should have the right level of sophistication to represent a stochastic population and therefore deal with multi-object estimation.

\subsection*{Q6: What is a ``multi-variate'' random experiment?}

\begin{figure}[h]
\centering
\def\svgwidth{0.5\columnwidth}
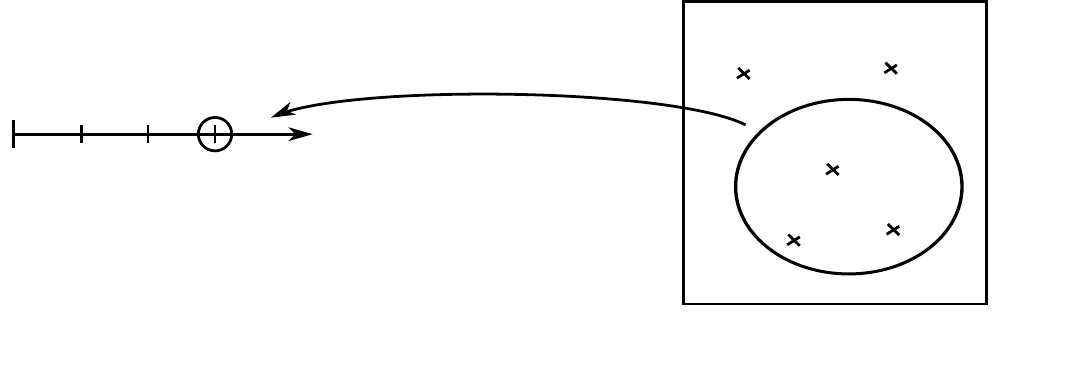
\caption{A counting function $N_{\cdot}(B)$ at $\vphi = \{x_1,\ldots,x_5\}$.}
\label{fig:countingFunction}
\end{figure}

A multi-variate random experiment consists of an unknown random number of random individuals represented by points, often referred to as a point process. There are several ways of defining a point process:

\begin{enumerate}
\itemName{Daley and Vere-Jones Vol.~I \cite{DaleyVere-Jones_2003}} A point process $\Phi$ can be characterised by:
\begin{enumerate}
\item A distribution $\{p_n\}_{n\geq 0}$, determining the total number of points in the population, with
\eqn{
\sum_{n\geq 0} p_n = 1.
}
\item For each integer $n \geq 1$, a probability measure $P^{(n)}_{\Phi}$ on $\mathcal{B}(\mathcal{X}^n)$, determining the joint distribution of the states of the points in the process, given that their total number is $n$.
\end{enumerate}

For every possible cardinality $n$, the point process $\Phi$ is then characterised on the space $\mathcal{X}^n$, i.e.\ a space of sequences or vectors rather than a space of sets. Because there is no natural order in most of the high dimensional spaces, e.g.\ $\mathcal{X} = \mathbb{R}^d$ with $d > 1$, the order of the points in the sequence has to be irrelevant. To ensure this is true, for any $n \in \mathbb{N}^*$, we can symmetrise the probability measure $P^{(n)}_{\Phi}$ and therefore define a measure $J^{(n)}_{\Phi}$, called Janossy measure, as follows
\eqn{
\label{eq:defJanossyMeasure}
J^{(n)}_{\Phi}(B_1 \times\ldots\times B_n) = p_n \sum_{\sigma \in S_n} P^{(n)}_{\Phi}(B_{\sigma(1)}\times\ldots\times B_{\sigma(n)}),
}
where $S_n$ is the permutation group on $n$ letters and $B_1 \times \ldots \times B_n \in \mathcal{B}(\mathcal{X}^n)$. The Janossy measure is a useful and insightful characterisation of a point process and has many properties.

Even though this definition of the point process $\Phi$ is sufficient to characterise it, we would like to isolate the randomness in another space, say $(\Omega,\mathcal{F},\mathbb{P})$, for the same reasons as for random variables.

\itemName{Daley and Vere-Jones Vol.~II \cite{DaleyVere-Jones_2008}} Point processes can be characterised by counting measures, which are defined as follows.

\begin{definition}[Counting measure] For a given set $\vphi$ of points in $\mathcal{X}$ such that $\vphi = \{x_1,\ldots,x_n\}$, the counting measure $N_{\vphi}$ is defined as
\eqnalign{
N_{\vphi}:\mathcal{B}(\mathcal{X}) & \to \mathbb{N} \\
\notag
B & \mapsto \sum_{i=1}^n \delta_{x_i}(B),
}
where $\delta_{x}$ is the Dirac measure at point $x$ such that $\delta_{x}(B) = 1$ if $x \in B$ and $0$ otherwise.
\end{definition}

Note that, as a measure, the Dirac measure $\delta_{x}$ can take a function as an argument, and be written $\delta_{x}(f) = f(x)$.

A point process can then be defined as a random measure:
\eqnalign{
N: (\Omega,\mathcal{F},\mathbb{P}) & \to (\mathcal{N}^{\#*}_{\mathcal{X}}, \mathcal{B}(\mathcal{N}^{\#*}_{\mathcal{X}}))\\
\notag
\omega & \mapsto N_{\omega},
}
where $\mathcal{N}^{\#*}_{\mathcal{X}}$ is the space of nice (simple boundedly finite) counting measure on $\mathcal{X}$.

\begin{itemize}
\item {\em Advantages}: inherit the properties and theorems related to random measures,
\item {\em Drawbacks}: not as intuitive as the first definition.
\end{itemize}

\itemName{Stoyan et al. \cite{Stoyan_D_1995}} Let $E$ be the space of sets of points in $\mathcal{X}$. The space $E$ can be expressed as
\eqn{
E = \phi \cup \bigcup_{n \geq 1} \mathcal{X}^{(n)},
}
where $\phi$ is the empty configuration and $\mathcal{X}^{(n)}$ is the space of sets of size $n$. We additionally assume that for any set $\vphi \in E\setminus \phi$, the following properties are verified:
\begin{itemize}
\item $\vphi$ is locally finite (each bounded subset of $\mathcal{X}$ must contain only a finite number of points of $\vphi$),
\item $\vphi$ is simple, i.e.
\eqn{\forall x_i,x_j \in \vphi,\quad x_i = x_j \implies i = j.}
\end{itemize}

A point process can then be defined as a measurable mapping
\eqn{
\label{eq:StoPointProcess}
\Phi : (\Omega,\mathcal{F},\mathbb{P}) \to (E,\mathcal{B}(E)),
}
and the associated counting function is, for any $B \in \mathcal{B}(\mathcal{X})$,
\eqnalign{
\label{eq:StoCountingFunc}
N_{\cdot}(B) : (E,\mathcal{B}(E)) & \to \mathbb{N}\\
\notag
\vphi & \mapsto N_{\vphi}(B).
}
An example of counting function is given in Figure \ref{fig:countingFunction}.
\end{enumerate}

Henceforth, the definition of Stoyan will be used since it is the closest to Finite Set Statistics (FiSSt).
Note that the point process $\Phi$, as defined in \eqref{eq:StoPointProcess}, and the counting function $N_{\cdot}(B)$, defined in \eqref{eq:StoCountingFunc}, can be composed to form the random variable $N_{\Phi}(B) : (\Omega,\mathcal{F},\mathbb{P}) \to \mathbb{N}$.

We now have objects of different natures:
\begin{itemize}
\item $N_{\Phi}\hphantom{(B)}$\quad: Random measure in $\mathcal{N}^{\#*}_{\mathcal{X}}$,
\item $N_{\Phi}(B)$\quad: Random variable in $\mathbb{N}$,
\item $N_{\vphi}\hphantom{(B)}$\quad: Measure on $\mathcal{B}(\mathcal{X})$.
\end{itemize}

Note that, since $N_{\Phi}$ and $N_{\vphi}$ are measures, we can consider measuring functions with them. For instance, for any measurable function $f$ on $\BorelMeasSpace{\mathcal{X}}$,
\subeqn{eq:countMeasOnFunc}{
N_{\Phi}(f) & = \int f(x) N_{\Phi}(\dd x) \\
& = \sum_{x \in \Phi} \delta_{x}(f) \\
& = \sum_{x \in \Phi} f(x)
}

As for random variables, we would like to define the expectation with respect to the point process $\Phi$.
\begin{enumerate}
\itemName{First try} One possibility is to transpose the definition for random variable, {\em mutatis mutandis}, and write
\subeqn{eq:badExpectGlobal}{
\mathbb{E}[\Phi(\cdot)] & = \int_{\Omega} \Phi(\omega) \mathbb{P}(\dd \omega)\\
\label{eq:badExpect}
& = \int_E \vphi \, P_{\Phi}(\dd \vphi),
}
where $P_{\Phi}$ is the pushforward probability measure of $\mathbb{P}$ into $E$. The integral \eqref{eq:badExpect} is not well defined in general since there is no natural way of summing sets together, for instance, the sum ``$\{x_1\} + \{x_2,x_3\}$'' is not properly defined. A different expression is then needed for the expectation.

\itemName{Second try} The other possibility is to reuse directly the expectation as defined before with the random variable $N_{\Phi}(B)$:
\subeqn{eq:expectPointProc}{
\mathbb{E}[N_{\Phi(\cdot)}(B)] & = \int_{\Omega} N_{\Phi(\omega)}(B) \mathbb{P}(\dd \omega)\\
& = \int_E N_{\vphi}(B) \, P_{\Phi}(\dd \vphi),\\
& = \mu^{(1)}_{\Phi}(B),
}
where $\mu^{(1)}_{\Phi}$ is the first order moment measure, or mean, or intensity measure, also denoted $\mu_{\Phi}$ for short. Unlike the mean of a random variable which is a real value, the mean of a point process is a function of space and gives the expected number of points in a given area $B$.
\end{enumerate}

For consistency, we can define a the equivalent of the Janossy measure $J_{\Phi}$ based on $P_{\Phi}$ such that
\eqn{
J_{\Phi}(B_1 \times \ldots \times B_n) = n! P_{\Phi}(B_1 \times \ldots \times B_n).
}
Note that, unlike the definition \eqref{eq:defJanossyMeasure}, there is no sum over permutations here since $P_{\Phi}$ is defined on a space of sets where there is no notion of order.

\begin{example}[Campbell theorem] The following useful result, called Campbell theorem, can be proved straightforwardly:
\subeqn{eq:campbellTheorem}{
\mathbb{E}\left[\,\sum_{x \in \Phi} f(x)\right] & = \int_E \left(\,\sum_{x \in \vphi} f(x)\right) P_{\Phi}(\dd \vphi)\\
& = \int_E \int_{\mathcal{X}} f(x) N_{\vphi}(\dd x) P_{\Phi}(\dd \vphi)\\
& = \int_{\mathcal{X}} f(x) \mu_{\Phi} (\dd x),
}
where relations \eqref{eq:countMeasOnFunc} and \eqref{eq:expectPointProc} have been used.
\end{example}

In Finite Set Statistics, another integral can be used in place of the Lebesgue integral over $E$, this integral is called the {\em set integral} and is defined as follows.
\begin{definition}[Set Integral]
Let $f$ be a multi-object density on $E$ and let $B_{\mathcal{X}}$ be a measurable subset in $\mathcal{B}(\mathcal{X})$. The set integral of $f$ on $B_{\mathcal{X}}$ is defined as
\eqn{
\int_{B_{\mathcal{X}}} f(\vphi) \delta\vphi = \sum_{n \geq 0} \dfrac{1}{n!} \int_{(B_{\mathcal{X}})^n} f(\{x_1,\ldots,x_n\}) \dd y_1 \ldots \dd y_n,
} 
\end{definition}
where $(B_{\mathcal{X}})^0 = \phi$ by convention.

\subsubsection*{Q6': What is the difference between the set integral and the Lebesgue integral on $E$?} 

Let $B_{\mathcal{X}} = B^{(1)}_{\mathcal{X}} \cup B^{(2)}_{\mathcal{X}}$ be a measurable subset of $\mathcal{X}$. The set integral of the multi-object density $f$ on $E$ is not additive in $B_{\mathcal{X}}$:
\eqn{
\int_{B_{\mathcal{X}}} f(\vphi) \delta \vphi  \neq \int_{B^{(1)}_{\mathcal{X}}} f(\vphi) \delta \vphi + \int_{B^{(2)}_{\mathcal{X}}} f(\vphi) \delta \vphi,
}
regardless of the fact that $B^{(1)}_{\mathcal{X}} \cap B^{(2)}_{\mathcal{X}} = \emptyset$ or not.

Let $B_E = B^{(1)}_E \cup B^{(2)}_E$ be a measurable subset of $E$ such that $B^{(1)}_E \cap B^{(2)}_E = \emptyset$. The set integral of the multi-object density $f$ on $E$ is additive in $B_E$:
\eqn{
\int_{B_E} \dfrac{1}{|\vphi|!} f(\vphi) \dd \vphi = \int_{B^{(1)}_E} \dfrac{1}{|\vphi|!} f(\vphi) \dd \vphi + \int_{B^{(2)}_E} \dfrac{1}{|\vphi|!} f(\vphi) \dd \vphi.
}

It is then a matter of choice: the Lebesgue integral is defined on the space of sets and might be complicated to handle while the set integral is defined on the individual space and is thus simpler. However, this simplicity comes at the cost of the loss of the additivity which is a considered as a usual property of the integral.

\chapter[Bayes' filter and PHD filter]{Bayes' filter and PHD filter}

In this lecture, the derivations of the Bayes' filter and of the Probability Hypothesis Density (PHD) filter will be given and explained. The objective is to see that, when using probability generating functionals along with relatively simple probability rules, we can arrive easily to the generating functional form of the prediction and the update. The similarities between these two steps will also be highlighted. The lecture begins with the introduction of the concept of probability generating functional starting from Fourier-like integral transforms.

\section{Generating functional}

Integral transform is popular concept in maths and engineering since it often allows for a simple and powerful representation of difficult concepts and the calculation of complicated quantities. Multi-object estimation is not an exception since the integral transforms called probability generating functionals (p.g.fl.s) can be used to greatly simplify some of the most challenging problems. The developments here follow 
\cite{DaleyVere-Jones_2008}:
\begin{itemize}
\item From random linear functional/random measures:\\
{\em Characteristic functional}:
\eqn{
F_{\Phi}[f] = \mathbb{E}\left[ \exp\left(\ii \int f\; \dd N_{\Phi}\right) \right].
}
\item Counting measures are non-negative, let $f = \ii g$:\\
{\em Laplace functional}:
\eqn{
\label{eq:LaplaceFunctional}
L_{\Phi}[g] = \mathbb{E}\left[ \exp\left(- \int g\; \dd N_{\Phi}\right) \right].
}
\item The Laplace functional is useful for calculating moments (see \cite{Stoyan_D_1995}, or Chapter \ref{chapt:HigherOrderMoments}). However, at this point, we are only interested in the simpler factorial moments and we can simplify the Laplace functional by setting $g = - \log h$:\\
{\em Probability generating functional (p.g.fl.)}:
\subeqn{eq:defGenFunc}{
G_{\Phi}[h] = L_{\Phi}[g] & = \mathbb{E}\left[ \exp\left( + \sum_{x \in \Phi} \log h(x) \right) \right] \\
\label{eq:defGenFunc_expect}
& = \mathbb{E}\left[ \; \prod_{x \in \Phi} h(x) \right] \\
& = \sum_{n \geq 0} \int_{\mathcal{X}^{(n)}} h(x_1)\ldots h(x_n) P_{\Phi}(\dd\left\{ x_1,\ldots,x_n \right\}).
}
\end{itemize}

\begin{example}\label{PgflSuperposition}Let $\Phi_1,\ldots,\Phi_N$ be independent point processes from $(\Omega,\mathcal{F},\mathbb{P})$ to $(E_X,\mathcal{B}(E_X))$ and let $\Phi = \Phi_1 \cup \ldots \cup \Phi_N$.

For a given outcome $\omega$,
\eqn{
\Phi(\omega) = \{x_{1,1},\ldots, x_{1,n_1}\} \cup \ldots \cup \{x_{N,1},\ldots, x_{N,n_N}\},
}
for some $n_1,\ldots,n_N \in \mathbb{N}$ and some $x_{i,j} \in \mathcal{X}$, $1\leq i \leq N$, $1\leq j \leq n_i$.
We can now rewrite the p.g.fl.\ $G_{\Phi}$ as a function of the p.g.fl.s $G_{\Phi_i}$, $1\leq i \leq N$:
\eqn{
\label{eq:PgflSuperposition}
G_{\Phi}[h] = \mathbb{E}\left[ \; \prod_{x \in \Phi} h(x) \right] = \mathbb{E}\left[ \; \prod_{i=1}^N \prod_{x \in \Phi_i} h(x) \right]
 = \prod_{i=1}^N G_{\Phi_i}[h].
}

This example shows that p.g.fl.s are useful to handle processes defined as a superposition of other independent processes.
\end{example}

We can still define Janossy measures, i.e.\ for any $B \in \mathcal{B}(E_X)$,
\subeqn{eq:JanossyMeasures}{
\label{eq:JanossyMeasures_1}
\int_B J_{\Phi}(\dd \{x_1,\ldots,x_n\}) & = \int_B n! P_{\Phi}(\dd\{x_1,\ldots,x_n\}) \\
\label{eq:JanossyMeasures_2}
& = \int_B f_{\Phi}(\{x_1,\ldots,x_n\}) \lambda(\dd\{x_1,\ldots,x_n\})\\
& = \int_B f_{\Phi}(\{x_1,\ldots,x_n\}) \dd x_1 \ldots \dd x_n \\
& = \int_B f_{\Phi}(X) \dd X,
}
where $\lambda$ is a reference measure on $E_X$ and where the Radon-Nikodym theorem has been used from \eqref{eq:JanossyMeasures_1} to \eqref{eq:JanossyMeasures_2} under the assumption that $P_{\Phi}$ is absolutely continuous with respect to $\lambda$ ($P_{\Phi} \ll \lambda$). The function $f_{\Phi}$ is called a multi-object density.

We can rewrite the p.g.fl.\ $G_{\Phi}$ of the point process $\Phi$ with the multi-object density $f_{\Phi}$ and get
\subeqn{eq:pgflWithJanossy}{
G_{\Phi}[h] & = \mathbb{E}\left[ \; \prod_{x \in \Phi} h(x) \right] \\
& = \sum_{n \geq 0} \int_{\mathcal{X}^{(n)}} \dfrac{1}{n!} h(x_1)\ldots h(x_n) f_{\Phi}(\{x_1,\ldots,x_n\}) \dd x_1 \ldots \dd x_n
}

\section{Multi-object estimation}

All the results in this section can be found either in \cite{Mahler_RPS_2007_3} or in \cite{Mahler_RPS_2003}, only the presentation and the notations differ.

\subsection{Prediction: Chapman-Kolmogorov equation}

The objective in this section is to relate the spaces $E_X$ and $E_Y$, defined as follows
\eqn{
E_X = \phi \cup \bigcup_{n \geq 1} \mathcal{X}^{(n)}, \qquad \qquad E_Y = \phi \cup \bigcup_{n \geq 1} \mathcal{Y}^{(n)}.
}

More precisely, we want to move to the space $E_X$ all the knowledge we have about the point process $\Phi_Y$ in $E_Y$, which is represented by the multi-object density $f_Y$.

\subsubsection{Bayes' filter prediction}

The most direct way of relating two point processes is to define a joint multi-object density $f_{X,Y}$ on $E_X \times E_Y$. It is useful to express $f_{X,Y}$ as a function of $f_Y$:
\eqn{
f_{X,Y}(X,Y) = M_{X|Y}(X|Y)f_Y(Y),
}
where $M_{X|Y}$ is a multi-object Markov transition. The next step consists of extracting from $f_{X,Y}$ all the information related to $\Phi_X$. More formally, the objective is to {\em marginalise} the joint multi-object density $f_{X,Y}$ as follows
\subeqn{eq:CKMargJoint}{
f_X(X) & = \int_{E_Y} \dfrac{1}{|Y|!} f_{X,Y} (X,Y) \dd Y\\
& = \int_{E_Y} \dfrac{1}{|Y|!} M_{X|Y}(X|Y)f_Y(Y) \dd Y.
}
This is the Chapman-Kolmogorov equation for point processes. The term $1/|Y|!$ in \eqref{eq:CKMargJoint} can be surprising at first sight, however, since $f_{X,Y}$ is not a probability density and since marginalisation applies only on probability densities, it is necessary to renormalise $f_{X,Y}$ before integrating it. According to \eqref{eq:JanossyMeasures_1}, the normalising factor of $f_{X,Y}(\cdot,Y)$ is $|Y|!$.

\subsubsection{P.g.fl.\ form of the Chapman-Kolmogorov equation}

We have seen that probability generating functionals (p.g.fl.s) are useful for dealing with point processes, and we would like to introduce them in \eqref{eq:CKMargJoint}. This can be done by multiplying on the left and right by $\frac{1}{|X|!}\left( \prod_{x\in X} h(x) \right)$ and integrating over $E_X$:
\xline{
\int_{E_X} \dfrac{1}{|X|!} \left( \prod_{x\in X} h(x) \right) f_X(X) \dd X = \\
\int_{E_X}\int_{E_Y} \dfrac{1}{|X|!|Y|!} \left( \prod_{x\in X} h(x) \right) M_{X|Y}(X|Y)f_Y(Y) \dd Y \dd X,
}
where the p.g.fl.s $G_X$ and $G_{X|Y}$ can be identified so that
\eqn{
\label{eq:CKpgflForm_1}
G_X[h] = \int_{E_Y} \dfrac{1}{|Y|!} G_{X|Y}[h|Y]f_Y(Y) \dd Y.
}

Given that $f_Y$ is known, the only unknown in \eqref{eq:CKpgflForm_1} is the p.g.fl.\ $G_{X|Y}$. Assuming that all the targets are independent of each other, that the birth process $\Phi_{\birth}$ is known and represented by $G_{\birth}$, and that there is no spawning:
\begin{itemize}
\item if $Y=\phi$, $G_{X|Y}[h|\phi] = G_{\birth}[h]$,
\item if $Y = \{y\}$, no birth,
\subeqn{eq:Gtarget}{
G_{X|Y}[h|\{y\}] & = \int_{E_X} \dfrac{1}{|X|!} \left( \prod_{x\in X} h(x) \right) M_{X|Y}(X|\{y\}) \dd X \\
& = M_{X|Y}(\phi|\{y\}) + \int_{\mathcal{X}} h(x) M_{X|Y}(\{x\}|\{y\}) \dd x \\
& = G_{\target}[h|y],
}
Since the multi-object Markov transition $M_{X|Y}$ is defined on $E_X \times E_Y$ and is a multi-object density on $E_X$, we can write
\eqn{
G_{\target}[1|y] = 1 = M_{X|Y}(\phi|\{y\}) + \int_{\mathcal{X}} M_{X|Y}(\{x\}|\{y\}) \dd x.
}
Let $\hat{M}_{X|Y}$ be a Markov transition on $\mathcal{X}\times\mathcal{Y}$ and $p_S$ a function on $\mathcal{Y}$ such that $M_{X|Y}(\{x\}|\{y\}) = p_S(y)\hat{M}_{X|Y}(x|y)$, then $M_{X|Y}(\phi|\{y\}) = 1 - p_S(y)$, and
\eqn{
G_{\target}[h|y] = 1 - p_S(y) + \int_{\mathcal{X}} h(x) p_S(y)\hat{M}_{X|Y}(x|y) \dd x.
}
\item if $Y = \{y_1,\ldots,y_n\}$, the point process $\Phi_X(Y)$, conditioned on $Y$, can be expressed as
\eqn{
\Phi_X(Y) = \Phi_{\birth} \cup \Phi_{X,1}(y_1) \cup \ldots \cup \Phi_{X,n}(y_n),
}
so that, using the result of Example \ref{PgflSuperposition},
\eqn{
\label{eq:CKfullConditionalPgfl}
G_{X|Y}[h|Y] = G_{\birth}[h]\prod_{i=1}^n G_{\target}[h|y_i].
}
\end{itemize}

We can now replace $G_{X|Y}[h|Y]$ in \eqref{eq:CKpgflForm_1} by its expression \eqref{eq:CKfullConditionalPgfl} and get
\eqn{
G_X[h] = G_{\birth}[h]\int_{E_Y} \dfrac{1}{|Y|!} \Big( \prod_{y \in Y} G_{\target}[h|y]\Big) f_Y(Y) \dd Y,
}
where we can identify the p.g.fl.\ $G_Y$ of the process $\Phi_Y$ and write
\eqn{
\label{eq:CKpgfl}
G_X[h] = G_{\birth}[h] G_Y[G_{\target}[h|\cdot]].
}

Equation \eqref{eq:CKpgfl} is the p.g.fl.\ form of the Chapman-Kolmogorov equation representing the prediction step.

\subsubsection{The PHD filter prediction}

The prediction equation of the PHD filter can be found straightforwardly by using \eqref{eq:CKpgfl} and given that
\eqn{
\mu_X(x) = \left.\delta G_X[h;\delta_x]\right|_{h=1}.
}
Indeed, by using the product rule for functional derivatives,
\eqn{
\label{eq:derivPHDpred}
\mu_X(x) = \left.\delta G_{\birth}[h;\delta_x]\right|_{h=1}G_Y\big[G_{\target}[1|\cdot]\big] + G_{\birth}[1] \left.\delta\left( G_Y\big[G_{X|Y}[h|\cdot]\big];\delta_x \right)\right|_{h=1}.
}
Since any p.g.fl.\ $G$ satisfy $G[1] = 1$, we can see that
\subeqn{eq:CKpgflAtOne}{
G_{\birth}[1] & = 1, \\
G_{\target}[1|\cdot] & = 1, \\
G_Y\big[G_{\target}[1|\cdot]\big] & = G_Y[1] = 1.
}

Since the birth process $\Phi_{\birth}$ is assumed to be known, its first moment density $\mu_{\birth}(x) =  \left.\delta G_{\birth}[h;\delta_x]\right|_{h=1}$ is also considered as known. The last unknown term in \eqref{eq:derivPHDpred} is the derivative of $G_Y[G_{X|Y}[h|\cdot]]$ in the direction $\delta_x$ at point $h=1$. If, following \eqref{eq:defGenFunc_expect}, we rewrite $G_Y$ as an expectation, then this term becomes
\subeqn{eq:derivPHDpred_2}{
\delta\big( G_Y\big[G_{X|Y}[h|\cdot]\big] & ;\delta_x \big)\big|_{h=1} \\
& = \mathbb{E}\Big[ \delta\Big( \prod_{y \in \Phi_Y} G_{\target}[h|y]; \delta_x\Big)  \Big]\Big|_{h=1}\\
& = \mathbb{E}\Big[ \sum_{y \in \Phi_Y} \delta G_{\target}[h|y; \delta_x] \prod_{\substack{y' \in \Phi_Y \\ y' \neq y}} G_{\target}[1|y'] \Big] \\
& = \mathbb{E}\Big[ \sum_{y \in \Phi_Y} p_S(y) \hat{M}_{X|Y}(x|y) \Big] \\
& = \int_{\mathcal{Y}} p_S(y) \hat{M}_{X|Y}(x|y) \mu_Y(y) \dd y,
}
where the last line is due to Campbell theorem \eqref{eq:campbellTheorem}.

Rewriting \eqref{eq:derivPHDpred}, the {\em predicted} first moment density $\mu_X$ can be expressed in terms of the {\em prior} first moment density $\mu_Y$ as
\eqn{
\mu_X(x) = \mu_{\birth}(x) + \int_{\mathcal{Y}} p_S(y) \hat{M}_{X|Y}(x|y) \mu_Y(y) \dd y.
}

This last result is the PHD filter prediction \cite{Mahler_RPS_2003}.

\subsubsection{Example of Bayes' filter prediction}
\label{sssec:predictionExample}

We will only study a special case of the Bayes' filter prediction since the general case is rather involved. Indeed, if we were to differentiate the conditional p.g.fl.\ $G_{X|Y}$, given in \eqref{eq:CKfullConditionalPgfl}, $m$ times in the directions $\delta_{x_1},\ldots,\delta_{x_m}$ at point $0$, we would find the multi-object Markov density $M_{X|Y}$ to be
\eqn{
\label{eq:fullMarkov}
M_{X|Y}(X|Y) = f_{\birth}(X) \prod_{i=1}^n (1-p_S(y_i)) \sum_{\theta} \prod_{i: \theta(i) > 0} \dfrac{p_S(y_i)\hat{M}_{X|Y}(x_{\theta(i)}|y_i)}
{(1-p_S(y_i))\mu_{\birth}(x_i)},
}
where $X = \{x_1,\ldots,x_m\}$, $Y = \{y_1,\ldots,y_n\}$, where the birth is Poisson:
\eqn{
f_{\birth}(X) = \exp\left(-\int \mu_{\birth}(x) \dd x\right)\mu_{\birth}(x_1)\ldots\mu_{\birth}(x_m),
}
and where $\theta$ is a function from $\{1,\ldots,n\}$ to $\{0,1,\ldots,m\}$ such that $\theta(i) = \theta(j) \neq 0$ implies that $i \neq j$. The function $\theta$ can be understood as an association function. This is surprising since we do not usual think about association in the prediction step. As we will see below, this association usually disappear because of the symmetry of the prior multi-object density $f_Y$.

We assume that there is no birth and no death so that, if there is any $i \in \{1,\ldots,n\}$ such that $\theta(i) = 0$, then $M_{X|Y}(X|Y) = 0$. The consequences are that $m=n$ and that the only mappings $\theta$ left are the permutations $\sigma$ in $S_n$. The multi-object Markov density $M_{X|Y}$ becomes
\eqn{
\label{eq:simpleMarkov}
M_{X|Y}(X|Y) = \sum_{\sigma \in S_n} \prod_{i=1}^n \hat{M}_{X|Y}(x_{\sigma(i)}|y_i).
}

Using \eqref{eq:CKMargJoint} and \eqref{eq:simpleMarkov}, the {\em predicted} multi-object density $f_X$ can be expressed as
\eqn{
f_X(\{x_1,\ldots,x_n\}) = \int_{\mathcal{Y}^{(n)}} \dfrac{1}{n!} \sum_{\sigma \in S_n} \Big( \prod_{i=1}^n \hat{M}_{X|Y}(x_{\sigma(i)}|y_i)\Big) f_Y(\{y_1,\ldots,y_n\}) \dd y_1 \ldots \dd y_n.
}

Because of the symmetry of $f_Y$, the $n!$ permutations in $S_n$ are all equal and the Bayes' filter prediction can be written
\eqn{
f_X(\{x_1,\ldots,x_n\}) = \int_{\mathcal{Y}^{(n)}} \hat{M}_{X|Y}(x_1|y_1)\ldots\hat{M}_{X|Y}(x_n|y_n) f_Y(\{y_1,\ldots,y_n\}) \dd y_1 \ldots \dd y_n.
}

Note that multi-object densities are ``stable'' under prediction, in other words, the operation of prediction preserves its structure. The prediction for probability densities would be
\eqn{
p_X(\{x_1,\ldots,x_n\}) = \dfrac{1}{n!}\int_{\mathcal{Y}^{(n)}} \hat{M}_{X|Y}(x_1|y_1)\ldots\hat{M}_{X|Y}(x_n|y_n) p_Y(\{y_1,\ldots,y_n\}) \dd y_1 \ldots \dd y_n,
}
where the additional $1/n!$ makes the expression less intuitive.

\subsection{Update: Bayes' theorem}

We are given some kind of information about our random experiment under the form of the realisation $Z = \{z_1,\ldots,z_m\}$ of a point process $\Phi_Z$ in $E_Z$. We now want to relate the spaces $E_X$ and $E_Z$ in order to improve our knowledge about the  random experiment in $E_X$ using the realisation $Z$ of $\Phi_Z$.

\subsubsection{The Bayes' filter update}

As for the Chapman-Kolmogorov equation, we define a joint multi-object density $f_{Z,X}$ on $E_Z \times E_X$ and express it as
\eqn{
\label{eq:jointMeasState}
f_{Z,X}(Z,X) = L_{Z|X}(Z|X)f_X(X),
}
where $L_{Z|X}$ is understood as a multi-object likelihood on $E_Z \times E_X$. If we rewrite \eqref{eq:CKMargJoint} but with $f_{Z,X}$, we obtain
\eqn{
\label{eq:BayesDenom}
f_Z(Z) = \int_{E_X} \dfrac{1}{|X|!} L_{Z|X}(Z|X)f_X(X) \dd X.
}

This equation tells us how likely the realisation $Z$ is, considering all the possible states of the point process $\Phi_X$. However, this is no longer the quantity of interest. Indeed, we would rather like to have an expression for the conditional multi-object density $f_{X|Z}$. This can be easily obtained by developing the left hand side of \eqref{eq:jointMeasState} as follows
\eqn{
f_{Z,X}(Z,X) = f_{X|Z}(X|Z)f_Z(Z) = L_{Z|X}(Z|X)f_X(X).
}

We can actually reuse \eqref{eq:BayesDenom} and write:
\eqn{
\label{eq:BayesUpdate}
f_{X|Z}(X|Z) = \dfrac{L_{Z|X}(Z|X)f_X(X)}{\int_{E_X} \frac{1}{|X'|!} L_{Z|X}(Z|X')f_X(X') \dd X'}.
}
This is Bayes' theorem for point processes. As before, we would like a p.g.fl.\ form for it.

\subsubsection{P.g.fl.\ form of Bayes' theorem}

Multiplying \eqref{eq:BayesUpdate} by $\frac{1}{|X|!}\left(\prod_{x\in X} h(x) \right)$ and integrating over $E_X$ gives
\eqn{
\label{eq:BayesPgflForm_1}
G_{X|Z}[h|Z] = \dfrac{\int_{E_X} \frac{1}{|X|!} \left(\prod_{x\in X} h(x) \right) L_{Z|X}(Z|X)f_X(X) \dd X}
{\int_{E_X} \frac{1}{|X'|!} L_{Z|X}(Z|X')f_X(X') \dd X'}.
}

The question is now: what is the expression of the likelihood $L_{Z|X}$? It can be expressed as a $m^{th}$-order derivative
\eqn{
\label{eq:BayesPgflForm_2}
L_{Z|X}(Z|X) = \left.\delta^m G_{Z|X}[g|X;\delta_{z_1},\ldots,\delta_{z_m}]\right|_{g = 0}.
}

This does not solve our problem, since we still have to find the conditional p.g.fl.\ $G_{Z|X}$. However, assuming that targets generate independent observations, with no more than one observation per target, that the clutter process $\Phi_{\clutter}$ is known and represented by its p.g.fl.\ $G_{\clutter}[g]$:
\begin{itemize}
\item if $X=\phi$, $G_{Z|X}[g] = G_{\clutter}[g]$,
\item if $X = \{x\}$, no clutter,
\subeqn{eq:Gobs}{
G_{Z|X}[g|\{x\}] & = \int_{E_Z} \dfrac{1}{|Z|!} \left( \prod_{z\in Z} g(z) \right) L_{Z|X}(Z|\{x\}) \dd Z \\
& = L_{Z|X}(\phi|\{x\}) + \int_{\mathcal{Z}} g(z) L_{Z|X}(\{z\}|\{x\}) \dd z \\
& = G_{\obs}[g|x],
}
Since the multi-object likelihood $L_{Z|X}$ is defined on $E_Z \times E_X$ and since $L_{Z|X}(\cdot|X)$ is a multi-object density on $E_Z$ for any $X \in E_X$, we can write
\eqn{
G_{\obs}[1|x] = 1 = L_{Z|X}(\phi|\{x\}) + \int_{\mathcal{Z}} L_{Z|X}(\{z\}|\{x\}) \dd x.
}
Let $\hat{L}_{Z|X}$ be a likelihood on $\mathcal{Z}\times\mathcal{X}$ and $p_D$ a function on $\mathcal{X}$ such that $L_{Z|X}(\{z\}|\{x\}) = p_D(x)\hat{L}_{Z|X}(z|x)$, then $L_{Z|X}(\phi|\{x\}) = 1 - p_D(x)$, and
\eqn{
G_{\obs}[g|x] = 1 - p_D(x) + \int_{\mathcal{Z}} g(z) p_D(x)\hat{L}_{Z|X}(z|x) \dd z.
}
\item if $X = \{x_1,\ldots,x_n\}$, the point process $\Phi_Z(X)$, conditioned on $X$, can be expressed as
\eqn{
\Phi_Z(X) = \Phi_{\clutter} \cup \Phi_{Z,1}(x_1) \cup \ldots \cup \Phi_{Z,n}(x	_n),
}
so that, using the result of Example \ref{PgflSuperposition},
\eqn{
\label{eq:BayesFullConditionalPgfl}
G_{Z|X}[g|X] = G_{\clutter}[g]\prod_{i=1}^n G_{\obs}[g|x_i].
}
\end{itemize}

There are striking similarities between the p.g.fl.\ of the likelihood \eqref{eq:BayesFullConditionalPgfl} and the p.g.fl.\ form of the multi-object Markov transition \eqref{eq:CKfullConditionalPgfl}. Actually, they have the same structure and we can identify the equivalent terms:
\begin{align*}
\Phi_{\clutter}& \lrarrow \Phi_{\birth} \\
G_{\obs}& \lrarrow G_{\target}\\
L_{X|Y}, \hat{L}_{X|Y}& \lrarrow M_{X|Y}, \hat{M}_{X|Y}\\
p_D& \lrarrow p_S.
\end{align*}

In order to find the expression of the conditional p.g.fl.\ $G_{X|Z}[h|Z]$ in terms of p.g.fl.\ only, we can replace $G_{Z|X}[g|X]$ in \eqref{eq:BayesPgflForm_2} and \eqref{eq:BayesPgflForm_1} by its expression \eqref{eq:BayesFullConditionalPgfl} and get
\eqn{
\label{eq:BayesPgfl}
G_{X|Z}[h|Z] = \dfrac{\left.\delta^m G_{Z,X}[g,h;\delta_{z_1},\ldots,\delta_{z_m}]\right|_{g=0}}
{\left.\delta^m G_{Z,X}[g,1;\delta_{z_1},\ldots,\delta_{z_m}]\right|_{g=0}}
}
with
\eqn{
G_{Z,X}[g,h] = G_{\clutter}[g] G_X[h\,G_{\obs}[g|\cdot]].
}

Equation \eqref{eq:BayesPgfl} is the p.g.fl.\ form of the Bayes' theorem representing the update step.

As in Section \ref{sssec:predictionExample}, one can differentiate \eqref{eq:BayesFullConditionalPgfl} to find that the multi-object likelihood $L_{Z|X}$ is
\eqn{
\label{eq:fullLikeli}
L_{Z|X}(Z|X) = f_{\clutter}(Z) \prod_{i=1}^n (1-p_D(x_i)) \sum_{\theta} \prod_{i: \theta(i) > 0} \dfrac{p_D(x_i)\hat{L}_{Z|X}(z_{\theta(i)}|x_i)}
{(1-p_D(x_i))\mu_{\clutter}(z_i)},
}
where $Z = \{z_1,\ldots,z_m\}$, $X = \{x_1,\ldots,x_n\}$, where the clutter is Poisson:
\eqn{
f_{\clutter}(X) = \exp\left(-\int \mu_{\clutter}(x) \dd x\right)\mu_{\clutter}(x_1)\ldots\mu_{\clutter}(x_m),
}
and where $\theta$ is a function from $\{1,\ldots,n\}$ to $\{0,1,\ldots,m\}$ as defined in Section \ref{sssec:predictionExample}.

\subsubsection{PHD filter update}

To find the first order moment $\mu_{X|Z}$, one has to calculate the following derivative
\eqn{
\label{eq:firstMomentAsADiff}
\mu_{X|Z}(x) = \left.\delta G_{X|Z}[h|Z;\delta_x]\right|_{h=1}.
}

For this first moment to be easy to compute and in closed form, we have to assume that the prior multi-object density $f_X$ is Poisson, so that its p.g.fl.\ $G_X$ is
\eqn{
\label{eq:PoissonPgfl}
G_X[h] = \exp\left(\int (h(x) - 1)\mu_X(x) \dd x\right).
}

The result of \eqref{eq:firstMomentAsADiff} is then the update equation of the PHD filter \cite{Mahler_RPS_2003}:
\begin{multline}
\mu_{X|Z}(x) = \big(1-p_D(x)\big)\mu_X(x) \\
+ \sum_{z \in Z} \dfrac{p_D(x)\hat{L}_{Z|X}(z|x)\mu_X(x)}
{\mu_{\clutter}(z) + \int p_D(x')\hat{L}_{Z|X}(z|x')\mu_X(x') \dd x'}.
\end{multline}

\chapter[Modelling interactions]{Modelling interactions: Khinchin processes \\ \vspace{1 mm} {\Large Jeremie Houssineau and Daniel Clark}}

\section{A reminder about p.g.fl.}

\begin{definition}Let $\Phi$ be a point process with points in $\mathcal{X}$, described by the multi-object density $f_{\Phi}$. The p.g.fl.\ of $\Phi$ is a functional $G_{\Phi}$ defined as
\eqn{
G_{\Phi}:\mathcal{U} \to \mathbb{R}^+,
}
where $\mathcal{U}$ is the class of Borel measurable functions from $\mathcal{X}$ to $[0,1]$ with $1-h$ vanishing outside some bounded set. The p.g.fl.\ $G_{\Phi}$ can be expressed as
\eqn{
G_{\Phi}[h] = f_{\Phi}(\phi) + \sum_{n \geq 1} \int \dfrac{1}{n!} h(x_1)\ldots h(x_n) f_{\Phi}(\{x_1,\ldots,x_n\})\dd x_1\ldots \dd x_n,
}
with $\phi$ the empty configuration.
\end{definition}

Note that:
\begin{itemize}
\item $G_{\Phi}[1] = 1$,
\item $G_{\Phi}[0] = f_{\Phi}(\phi)$,
\item The multi-object density $f_{\Phi}$ is recovered, for any $n \in \mathbb{N}$ and any set $\{x_1,\ldots,x_n\} \in E_X$, through the following differential:
\eqn{
f_{\Phi}(\{x_1,\ldots,x_n\}) = \left. \delta^n G_{\Phi}[h; \delta_{x_1},\ldots,\delta_{x_n}]\right|_{h=0}.
}
\item The $n^{th}$-order factorial moment density $\alpha^{(n)}_{\Phi}$ is recovered, for any $n \in \mathbb{N}$ and any $x_1,\ldots,x_n \in \mathcal{X}$, through the following differential:
\eqn{
\label{eq:factorialMoment}
\alpha^{(n)}_{\Phi}(x_1,\ldots,x_n) = \left. \delta^n G_{\Phi}[h; \delta_{x_1},\ldots,\delta_{x_n}]\right|_{h=1}.
}
\end{itemize}

%
%

\begin{remark}[Interactions and correlations]
The meanings of the terms ``interaction'' and ``correlation'' are different. The former describes the action of one entity on another entity while the latter shows how dependent two random variables are. However, in the context of multi-object Bayesian estimation, correlations allow for interactions to be represented. Indeed, if two points $y_1$ and $y_2$ are correlated and represented by the multi-object density $f_Y$, one has to predict them jointly using a multi-object Markov density:
\eqn{
M_{X|Y}(x_1,x_2|y_1,y_2) f_Y(y_1, y_2).
}

Since $M_{X|Y}$ jointly propagates two points, e.g.\ in time, it can also change the state of one point considering the state of the other point. The concepts are therefore very close.
\end{remark}

The objective in the next sections is to build an example of interacting point process based on the example of the Poisson point process.

\section{A special case of interaction: no interaction}

In this section, we consider a well-known example of ``non-interacting'' point process: the Poisson point process.

Let $\Phi$ be a Poisson point process with intensity $\lambda s(x)$, with $\lambda \in \mathbb{R}^+$ and $s$ a probability density, and p.g.fl.
\eqn{
G_{\Phi}[h] = \exp\left( -\lambda + \lambda \int h(x)s(x) \dd x \right).
}

To emphasise the nature of this p.g.fl.\, the intensity can be defined as $K_{\Phi}^{(1)}(x) = \lambda s(x)$ and the p.g.fl.\ $G_{\Phi}$ can be rewritten as
\eqnalign{
G_{\Phi}[h] & = \exp\left( -\int K_{\Phi}^{(1)}(x) \dd x + \int h(x) K_{\Phi}^{(1)}(x) \dd x \right) \\
\label{eq:PoissonAsQuotient}
& = \dfrac{\exp\left(\int h(x) K_{\Phi}^{(1)}(x) \dd x \right)}{\exp\left( \int K_{\Phi}^{(1)}(x) \dd x \right)}.
}

The numerator of \eqref{eq:PoissonAsQuotient} can be understood as a compound generator since it duplicates the intensity $K_{\Phi}^{(1)}$ as many times as there are points in the process, for instance
\eqn{
f_{\Phi}(\{x_1,\ldots,x_n\}) = \dfrac{K_{\Phi}^{(1)}(x_1)\ldots K_{\Phi}^{(1)}(x_n)}{\exp\left( \int K_{\Phi}^{(1)}(x) \dd x \right)},
}
while the denominator is a normalising factor, making sure that $G_{\Phi}[1] = 1$.

Even though \eqref{eq:PoissonAsQuotient} is a insightful expression of the p.g.fl.\ of $\Phi$, it is useful to simplify it by setting
\eqn{
K_{\Phi}^{(0)} = \int K_{\Phi}^{(1)}(x) \dd x,
}
so that the p.g.fl.\ $P_{\Phi}$ can be simplified to
\eqn{
G_{\Phi}[h] = \exp\left( -K_{\Phi}^{(0)} + \int h(x) K_{\Phi}^{(1)}(x) \dd x  \right).
}

\section[Pairwise interaction]{Another special case of interaction: pairwise interaction}

Equation \eqref{eq:PoissonAsQuotient} is a useful way of defining a point process for the reasons given in the previous section. When trying to generate a point process with pairwise interactions, the following question can arise:
\begin{quote}
What if we ``replace the compound $K_{\Phi}^{(1)}(x)$ by another compound $K_{\Phi}^{(2)}(x_1,x_2)$ representing pairwise-correlated points''?
\end{quote}

Given that $K_{\Phi}^{(2)}$ is symmetric, we can define a p.g.fl.\ $G_{\Phi'}$ such that
\subeqn{eq:pairwisePgfl}{
G_{\Phi'}[h] & = \dfrac{\exp\left( \int h(x_1)h(x_2) K_{\Phi'}^{(2)}(x_1,x_2) \dd x_1 \dd x_2 \right)}
{\exp\left( \int K_{\Phi'}^{(2)}(x_1,x_2) \dd x_1 \dd x_2 \right) }\\
& = \exp\left( - K_{\Phi'}^{(0)} + \int h(x_1)h(x_2) K_{\Phi'}^{(2)}(x_1,x_2) \dd x_1 \dd x_2 \right),
}
where $K_{\Phi'}^{(0)} = \int K_{\Phi'}^{(2)}(x_1,x_2) \dd x_1 \dd x_2 $.

How to check that the p.g.fl.\ $G_{\Phi'}$ does represent the point process we are looking for? One possibility is to compute the multi-object density $f_{\Phi}$ and check it is meaningful.
\begin{enumerate}
\itemName{First try} Let's differentiate $G_{\Phi'}$ one time, and then two times and then try to see what form it takes in the general case
\eqn{
f_{\Phi} (\{y_1,\ldots,y_n\}) = \left.\delta^n G_{\Phi'}[h; \delta_{y_1},\ldots,\delta_{y_n}]\right|_{h=0}, \quad y_1,\ldots,y_n \in \mathcal{X},
}
The first derivative can be computed easily:
\subeqn{eq:pairwise1order}{
\delta G_{\Phi'}[h;\delta_{y_1}] & = \left( \int h(x_2) K_{\Phi'}^{(2)}(y_1,x_2) \dd x_2 + \int h(x_1) K_{\Phi'}^{(2)}(x_1,y_1) \dd x_1 \right) G_{\Phi'}[h]\\
& = \left(2\int h(x_2) K_{\Phi'}^{(2)}(y_1,x_2) \dd x_2 \right)G_{\Phi'}[h],
}
where the symmetry of $K_{\Phi'}^{(2)}$ has been used. This is no more complicated than the calculations required for the derivation of the PHD filter update since it only consists of differentiating a exponential. 

We now want to compute the second derivative:
\begin{multline}
\delta^2 G_{\Phi'}[h;\delta_{y_1},\delta_{y_2}] = 2 K_{\Phi'}^{(2)}(y_1,y_2) G_{\Phi'}[h]\\
+ \left(2\int h(x_2) K_{\Phi'}^{(2)}(y_1,x_2) \dd x_2 \right)\left(2\int h(x_2) K_{\Phi'}^{(2)}(y_2,x_2) \dd x_2 \right)G_{\Phi'}[h].
\end{multline}

It is getting a bit more complicated and we see that if we keep differentiating, the calculations will be more and more difficult.

\itemName{Second try} For the sake of compactness, we set 
\eqn{
A^{(2)}_{\Phi'}[h] = -K_{\Phi'}^{(0)} + \int h(x_1)h(x_2) K_{\Phi'}^{(2)}(x_1,x_2) \dd x_1 \dd x_2.
}

The p.g.fl.\ $G_{\Phi'}$ can be expressed as a composition
\eqn{
G_{\Phi'}[h] = \exp\left(A^{(2)}_{\Phi'}[h]\right).
}

The rule for differentiating $n$ times a composition is called Fa\`a di Bruno's formula \cite{Clark2013} and can be written
\begin{multline}
\label{eq:PairwiseFaa}
\delta^n G_{\Phi'}[h; \delta_{y_1},\ldots,\delta_{y_n}] = \\
\sum_{\pi \in \Pi(y_1,\ldots,y_n)} \delta^{|\pi|} \exp\left( A^{(2)}_{\Phi'}[h]; \delta^{|\omega|} A^{(2)}_{\Phi'}[h; \delta_{y}:y \in \omega]:\omega \in \pi \right)
\end{multline}

We can easily compute the different terms in \eqref{eq:PairwiseFaa} when $h = 0$:
\begin{itemize}
\item $A^{(2)}_{\Phi'}[0] = \exp\left( -K_{\Phi'}^{(0)} \right)$,
\item $\left.\delta A^{(2)}_{\Phi'}[h; \delta_{y_1}]\right|_{h=0} = \left( \int h(x_2) K_{\Phi'}^{(2)}(y_1,x_2) \dd x_2 \right)_{h=0} = 0$,
\item $\left.\delta^2 A^{(2)}_{\Phi'}[h;\delta_{y_1},\delta_{y_2}]\right|_{h=0} = 2 K_{\Phi'}^{(2)}(y_1,y_2)$,
\item $\left.\delta^n A^{(2)}_{\Phi'}[h;\delta_{y_1},\ldots,\delta_{y_n}]\right|_{h=0} = 0$, for any $n > 2$.
\end{itemize}

We can conclude that, if $n$ is even,
\eqn{
\label{eq:modOfPairwiseInter}
f_{\Phi}(\{y_1,\ldots,y_n\}) = \exp\left( -K_{\Phi'}^{(0)} \right) \sum_{\pi \in \Pi_2(y_1,\ldots,y_n)} \prod_{\omega \in \pi} 2 K_{\Phi'}^{(2)}(y_1,y_2),
}
where $\Pi_2(y_1,\ldots,y_n)$ is the set of binary partitions of $\{y_1,\ldots,y_n\}$. Equation \eqref{eq:modOfPairwiseInter} can be interpreted easily, the sum over binary partitions represent all the possible ways of defining couples out of $n$ points and for each of these couples, the probability of its existence and the associated probability distribution is given by $K_{\Phi'}^{(2)}$. 

However, if $n$ is odd, $f_{\Phi}(\{y_1,\ldots,y_n\}) = 0$. The conclusion is that $\Phi'$ is a point process containing only pairs of interacting points, so it makes sense that the probability for such a process to contain an odd number of points is $0$.

This is very constraining, we would like both independent {\em and} pairwise-correlated points. The solution is to superpose the two point processes that have been introduced before, the Poisson point process $\Phi$ and the pairwise-correlated point process $\Phi'$ to form a new point process $\Phi_{\text{GP}}$. This point process has already been defined before and bear the name {\em Gauss-Poisson} point process (hence justifying the subscript GP). It is straightforward to define the p.g.fl.\ $G_{\text{GP}}$ of the point process $\Phi_{\text{GP}}$ since it is a superposition: 
\eqn{
G_{\text{GP}}[h] = G_{\Phi}[h] G_{\Phi'}[h].
}

\begin{remark}[How to represent a Gauss-Poisson point process?] We know that the Poisson point process can be represented by its first moment density. Can the Gauss-Poisson point process be represented by moments or factorial moments as well?
\begin{itemize}
\item First moment density:
\subeqn{eq:firstMomentGP}{
\mu^{(1)}_{\text{GP}}(y) = \alpha^{(1)}_{\text{GP}}(y) & = \left.\delta G_{\Phi'}[h;\delta_{y_1}]\right|_{h=1}\\
& = K^{(1)}_{\text{GP}}(y) + 2 \int K^{(2)}_{\text{GP}}(y,x) \dd x.
}
The first moment density $\mu^{(1)}_{\text{GP}}$ is not sufficient to characterise a Gauss-Poisson point process since there are two unknown quantities $K^{(1)}_{\text{GP}}$ and $K^{(2)}_{\text{GP}}$ and a single equation.

\item Second factorial moment density:
\begin{align}
\notag
& \alpha^{(2)}_{\text{GP}}(y_1,y_2) = 2 K^{(2)}_{\text{GP}}(y_1,y_2) \\
& + \left( K^{(1)}_{\text{GP}}(y_1) + 2 \int K^{(2)}_{\text{GP}}(y_1,x) \dd x \right)\left( K^{(1)}_{\text{GP}}(y_2) + 2 \int K^{(2)}_{\text{GP}}(y_2,x) \dd x \right)
\end{align}
which can be rewritten in terms of the first moment density $\mu^{(1)}_{\text{GP}}$ as
\eqn{
\alpha^{(2)}_{\text{GP}}(y_1,y_2) = 2 K^{(2)}_{\text{GP}}(y_1,y_2) + \mu^{(1)}_{\text{GP}}(y_1)\mu^{(1)}_{\text{GP}}(y_2).
}
The second factorial moment density $\alpha^{(2)}_{\text{GP}}$ has a very specific form which makes natural the introduction of the covariance for point processes:
\subeqn{eq:covGP}{
\mathrm{Cov}_{\text{GP}}(y_1,y_2) & = \alpha^{(2)}_{\text{GP}}(y_1,y_2) - \mu^{(1)}_{\text{GP}}(y_1)\mu^{(1)}_{\text{GP}}(y_2)\\
& = 2 K^{(2)}_{\text{GP}}(y_1,y_2).
}
It is now clear that the covariance is directly related to the correlations in a point process. As a conclusion, a Gauss-Poisson point process can be represented by the pair mean/covariance $(\mu^{(1)}_{\text{GP}},\mathrm{Cov}_{\text{GP}})$.
\end{itemize}
\end{remark}

The idea illustrated in this section can be generalised to correlations of any order as described in the next section.

\section{The Khinchin point process}

\begin{definition}[Khinchin point process]
Let $\{K_{\Psi}^{(n)}\}_{n \geq 0}$ be a family of symmetric densities and let $\Psi$ be a point process defined through its p.g.fl.\ $G_{\Psi}$ expressed as
\eqn{
G_{\Psi}[h] = \exp\Big( - K_{\Phi'}^{(0)} + \sum_{n \geq 1} \int h(x_1)\ldots h(x_n) K_{\Phi'}^{(n)}(x_1,\ldots,x_n) \dd x_1 \ldots \dd x_n \Big),
}
where
\eqn{
K_{\Phi'}^{(0)} = \sum_{n \geq 1} \int K_{\Phi'}^{(n)}(x_1,\ldots,x_n) \dd x_1 \ldots \dd x_n \Big).
}
The point process $\Psi$ is referred to as Khinchin point process.
\end{definition}

To find the multi-object moment density of a Khinchin process $\Psi$ given its p.g.fl.\ $G_{\Psi}$, one can use Fa\`a di Bruno's formula the compute the $n^{th}$-order differential. The result of such an operation is
\eqn{
f_{\Psi}(\{x_1,\ldots,x_n\}) = \exp\left( -K^{(0)}_{\Psi} \right) \sum_{\pi \in \Pi(x_1,\ldots,x_n)} \prod_{\omega \in \pi} |\omega|! K^{(|\omega|)}_{\Psi}(\omega),
}
which can be interpreted in a way similar to \eqref{eq:modOfPairwiseInter}.

Note that, as a generalisation of the Poisson process, the Khinchin process generates i.i.d.\ compounds for a given cardinality, see Figure \ref{fig:KhinchinExample}.

We now have a basis for conducting Bayesian estimation with correlated point processes. Moreover, the p.g.fl.\ of a Khinchin process is based on the exponential which makes it easy to manipulate when equipped with the chain rule for higher-order differentials: Fa\`a di Bruno's formula.

\begin{figure}[h]
\centering
\def\svgwidth{0.5\columnwidth}
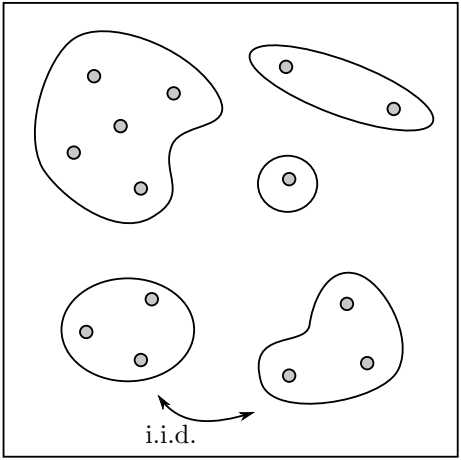
\caption{A realisation of a Khinchin process, with correlated points}
\label{fig:KhinchinExample}
\end{figure}

\end{enumerate}

\chapter[Higher-order moments]{Higher-order moments for point processes\\
\vspace{1 mm} {\Large Emmanuel Delande}}
\label{chapt:HigherOrderMoments}

This lecture exposes the concept of higher-order moments for point processes, and provides the tools for their practical derivation using the functional derivatives introduced in earlier lectures. 
The construction of the Probability Hypothesis Density (PHD) filter with variance in target number, an example of application of higher-order moments for target detection and tracking problems, is introduced. A more detailed construction is given in \cite{Delande_E_2013}, and a similar result for the Cardinalized Probability Hypothesis Density (CPHD) filter is established in \cite{Delande_E_2013_2}.
  
    \section{Point processes and counting measures}
      \begin{minipage}[t]{0.50\textwidth}
	\vspace{1cm}
	Moments are statistics describing \textit{random variables}. In our case, we need to build an integer-valued random variable providing a \textit{local} description of the target population.\newline
      \end{minipage}
      \begin{minipage}[t]{0.50\textwidth}
	\begin{figure}[H]
	  \centering
	  \includegraphics[width=0.65\textwidth]{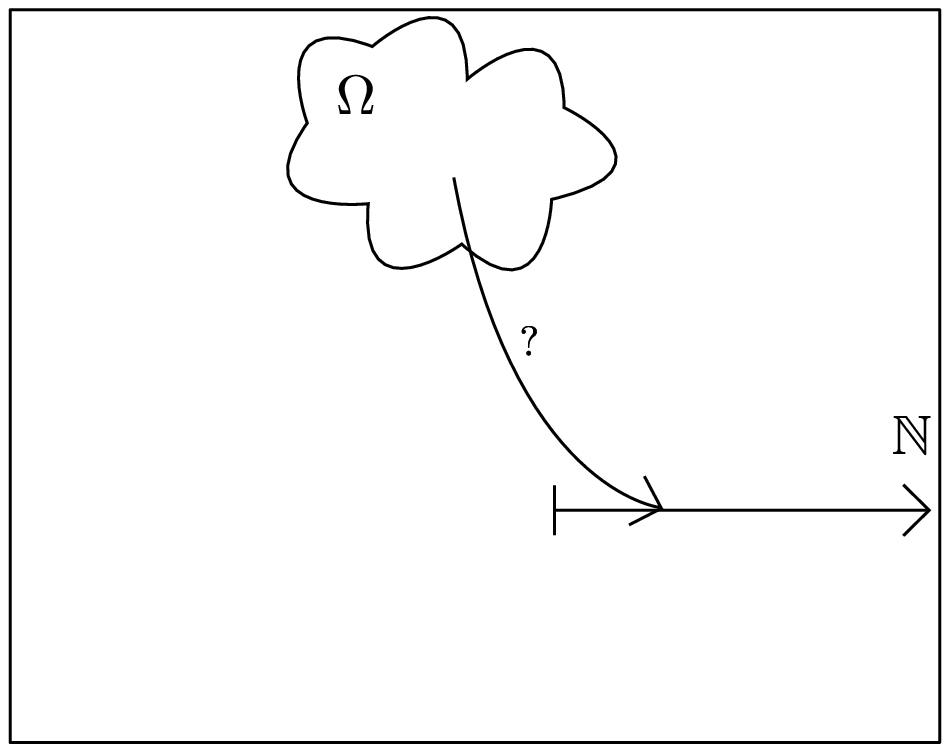}
	\end{figure}
      \end{minipage}
      
      \begin{minipage}[t]{0.50\textwidth}
	Point processes are \textit{random processes} whose realizations are set of points in the target state space $\Xcal$. A point process is \textit{not} a real (or complex) valued random variable and moments cannot be directly defined from point processes -- the expression $\Exp[\Phi]$ has no mathematical sense, since realizations can be sets of different sizes for which no sum operator is defined, see \eqref{eq:badExpectGlobal}.\newline
      \end{minipage}
      \begin{minipage}[t]{0.50\textwidth}
	\begin{figure}[H]
	  \centering
	  \includegraphics[width=0.65\textwidth]{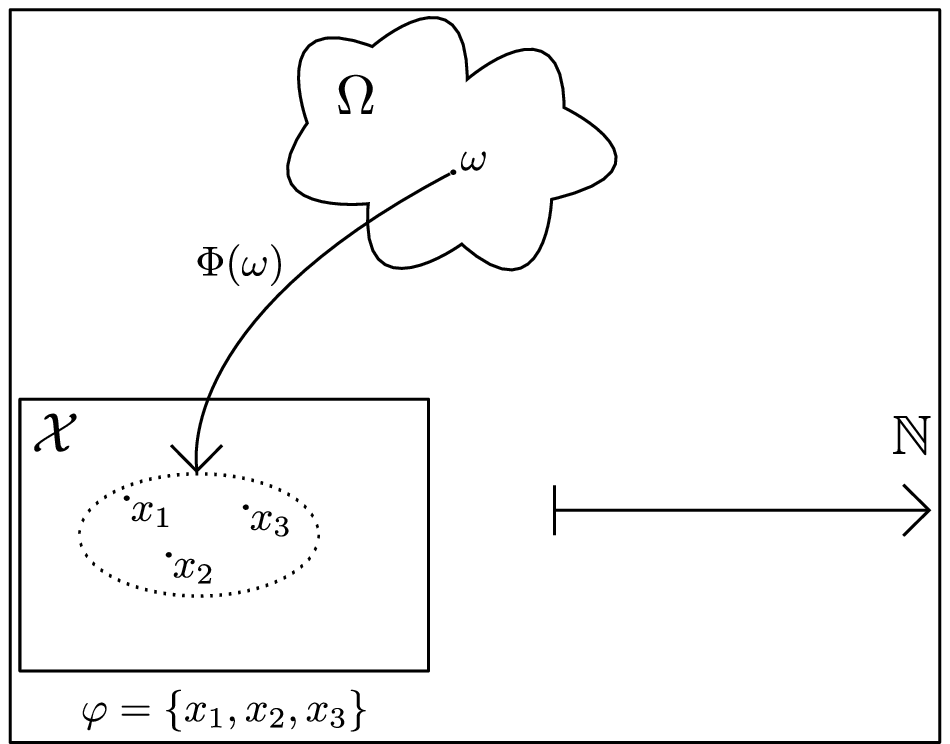}
	\end{figure}
      \end{minipage}
      
      \begin{minipage}[t]{0.50\textwidth}
	\vspace{1cm}
	Let us fix a region $B \in \Bcal(\Xcal)$\footnotemark[1]. One can map any realization $\varphi$ of the point process to the number of elements in $\varphi$ belonging to $B$, that is:
	\begin{equation} \label{EqCountingMeasure}
	  N_{\varphi}(B) = |\varphi \cap B|.
	\end{equation}   
      \end{minipage}
      \footnotetext[1]{$\Bcal(\Xcal)$, introduced in the previous lectures, see Chapter \ref{chapt:MeasureAndProba}, is the Borel $\sigma$-algebra on the target space $\Xcal$.}
      \begin{minipage}[t]{0.50\textwidth}
	\begin{figure}[H]
	  \centering
	  \includegraphics[width=0.65\textwidth]{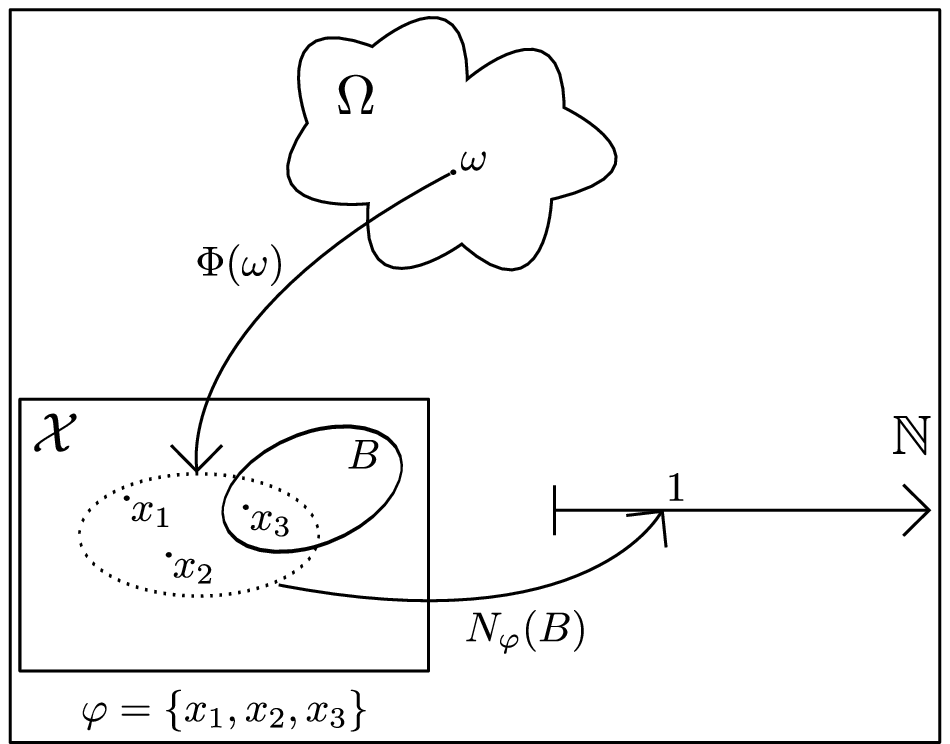}
	\end{figure}
      \end{minipage}
      
      \begin{minipage}[t]{0.50\textwidth}
	\vspace{0.5cm}
	If we compose the point process with the mapping defined above in \eqref{EqCountingMeasure}, we get the \textit{integer-valued random variable}
	\begin{equation} \label{EqRandomVariable}
	  N_{\Phi(\cdot)}(B) = |\Phi(\cdot) \cap B|,
	\end{equation}
	which provides a stochastic description of the number of targets \textit{inside $B$} acc. to the point process.
      \end{minipage}
      \begin{minipage}[t]{0.50\textwidth}
	\begin{figure}[H]
	  \centering
	  \includegraphics[width=0.65\textwidth]{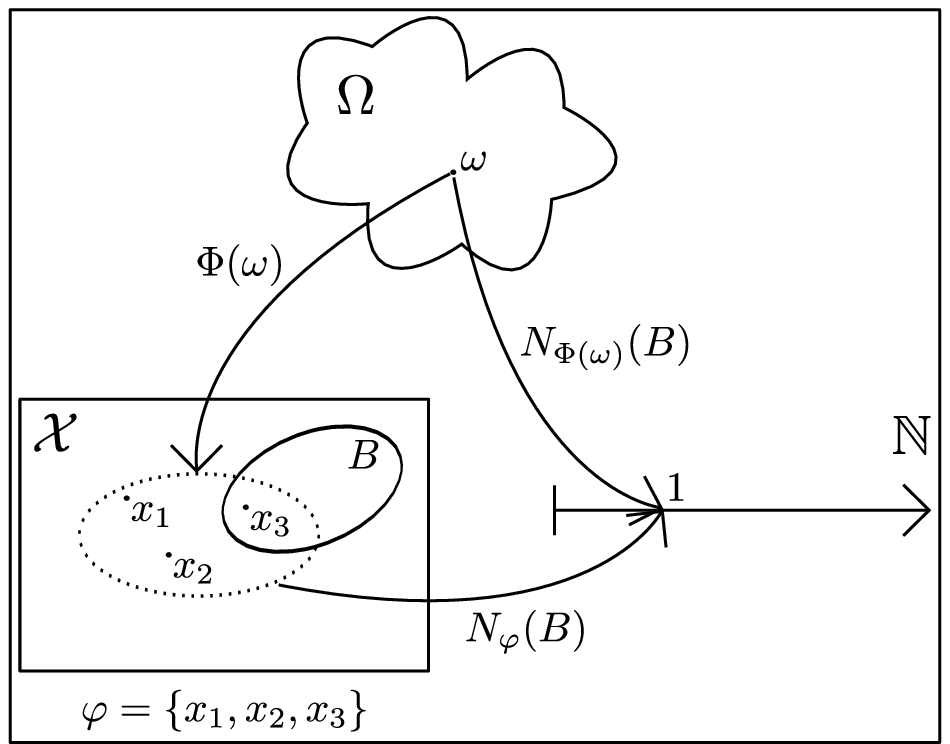}
	\end{figure}
      \end{minipage}
      
      \vspace{1cm}
      We have now built an integer-valued random variable $N_{\Phi(\cdot)}(B)$ for any region $B \in \Bcal(\Xcal)$, and we can now define its moments like for any real-valued random variable. First all, let us have a look at the nature of the different mathematical objects involved in the construction we have just illustrated.
      \begin{enumerate}
	\item $N_{\Phi(\omega)}(B)$ is the number of elements of the \textit{fixed} set $\Phi(\omega)$ belonging to the \textit{fixed} region $B$; it is an \textit{integer};
	\item $N_{\Phi(\cdot)}(B)$ maps an outcome $\omega \in \Omega$ to the number of elements of the set $\Phi(\omega)$ belonging to the \textit{fixed} region $B$; it is an \textit{integer-valued random variable};
	\item $N_{\Phi(\omega)}(\cdot)$ maps a region $B \in \Bcal(\Xcal)$ to the number of elements from the \textit{fixed} set $\Phi(\omega)$ that it contains; it is an \textit{integer-valued measure} called a \textit{counting measure};
	\item $N_{\Phi(\cdot)}(\cdot)$ maps an outcome $\omega \in \Omega$ to the counting measure $N_{\Phi(\omega)}(\cdot)$; it is a \textit{integer-valued random measure}.
      \end{enumerate}
      
    \section{Moment measures of point processes} \label{SecPointProcessMeasure}
      Since $N_{\Phi(\cdot)}(B)$ is a random variable for any region $B \in \Bcal(\Xcal)$, the joint expectation of any number of such quantities is well-defined. The \textit{$n$-th order moment measure} of the point process $\Phi$ is defined by the joint expectation
      \begin{multline} \label{EqMomentMeasure}
	\forall B_1,B_2,\dots,B_n \in \Bcal(\Xcal),
	\\
	\mu^{(n)}_{\Phi}(B_1,B_2,\dots,B_n) = \Exp[N_{\Phi(\cdot)}(B_1)N_{\Phi(\cdot)}(B_2) \cdots N_{\Phi(\cdot)}(B_n)].
      \end{multline}
      Note that by construction $\mu^{(n)}_{\Phi}$ is a function on the $\sigma$-algebra $\Bcal(\Xcal)^n$; in addition, it can be shown that it is a \textit{measure} on $\Bcal(\Xcal)^n$, hence its name. Note also that if $B_1 = \dots = B_n = B$ then we have
      \begin{equation} \label{EqMomentMeasureAndMoment}
	\mu^{(n)}_{\Phi}(B,B,\dots,B) = \Exp[\left(N_{\Phi(\cdot)}(B)\right)^n].
      \end{equation}
      That is, the real number $\mu^{(n)}_{\Phi}(B,B,\dots,B)$ is the $n$-th order moment of the random variable $N_{\Phi(\cdot)}(B)$.\newline
      \\
      We have now established that, for any region $B \in \Bcal(\Xcal)$:
      \begin{enumerate}
	\item $N_{\Phi(\cdot)}(B)$ is a random variable giving a statistical description of the number of targets within $B$;
	\item The $n$-th order moment of $N_{\Phi(\cdot)}(B)$ is given by $\mu^{(n)}_{\Phi}(B,B,\dots,B)$.
      \end{enumerate}
      While the knowledge of the random variable $N_{\Phi(\cdot)}(B)$ would provide a \textit{full} information of ``what is going on inside $B$'', it is usually unavailable in practical problems in which one must focus on the extraction of a few meaningful moments providing an \textit{approximate} information.\newline
      \\
      We shall focus from now on on the first two moment measures $\mu^{(1)}_{\Phi}$, $\mu^{(2)}_{\Phi}$, and the \textit{centered second moment} or \textit{variance} defined by:
      \begin{equation} \label{EqVariance}
	\forall B \in \Bcal(\Xcal),~var_{\Phi}(B) =  \mu^{(2)}_{\Phi}(B, B) - \left[\mu^{(1)}_{\Phi}(B)\right]^2.
      \end{equation}
      Note that by definition the variance is a function on the $\sigma$-algebra $\Bcal(\Xcal)$, but in general it is \textit{not} a measure on $\Bcal(\Xcal)$. This point is illustrated on a simple example later in the lecture slides.\newline
      \\
      Just as for any random variable, the statistics $(\mu^{(1)}_{\Phi}(\cdot), var_{\Phi}(\cdot))$ have the following interpretation:
      \begin{itemize}
	\item $\mu^{(1)}_{\Phi}(B)$ is the mean value of the random variable $N_{\Phi(\cdot)}(B)$;
	\item $var_{\Phi}(B)$ quantifies the spread of the random variable $N_{\Phi(\cdot)}(B)$ around its mean value.
      \end{itemize}
      In detection and tracking problems, these statistics allow us to estimate the average number of target in any desired region $B$, and provide an associated uncertainty.
      
    \section{Exploiting the Laplace functional} \label{SecLaplace}
      The first (non factorial) moment measure $\mu^{(1)}_{\Phi}$ is identical (see exercise \ref{SubsecLinearPGFlLaplace}) to the first \textit{factorial} moment measure  $\alpha^{(1)}_{\Phi}$ presented in a previous lecture, see e.g.\ \eqref{eq:factorialMoment}. Thus, $\mu^{(1)}_{\Phi}(B)$ can be retrieved from the first functional derivative of the PGFl $G_{\Phi}$:
      \begin{equation} \label{EqFirstMomentComputation}
	\forall B \in \Bcal(\Xcal),~\mu^{(1)}_{\Phi}(B) = \left.\delta G_{\Phi}(h; 1_B)\right|_{h = 1}
      \end{equation}  
      What about the second moment measure $\mu^{(2)}_{\Phi}$? Starting from the definition \eqref{EqMomentMeasure} we can write:
      \begin{subequations} \label{EqMomentMeasureDecomp1}
	\begin{align}
	  \mu^{(2)}_{\Phi}(B_1, B_2) &= \Exp[N_{\Phi}(B_1)N_{\Phi}(B_2)] \label{EqMomentMeasureDecomp1a} 
	  \\
	  &= \Exp\left[\sum_{x, x' \in \Phi}1_{B_1}(x)1_{B_2}(x')\right] \label{EqMomentMeasureDecomp1b} 
	  \\
	  &= \underbrace{\Exp\left[~\ds{\sideset{}{^{\neq}}\sum_{x, x' \in \Phi}1_{B_1}(x)1_{B_2}(x')}\right]}_{= \alpha^{(2)}_{\Phi}(B_1, B_2)} + \underbrace{\Exp\left[\sum_{x \in \Phi}1_{B_1}(x)1_{B_2}(x)\right]}_{= \mu^{(1)}_{\Phi}(B_1 \cap B_2)} \label{EqMomentMeasureDecomp1c} 
	\end{align}
      \end{subequations}
      That is, the second moment measure $\mu^{(2)}_{\Phi}$ can be retrieved from the second factorial moment measure $\alpha^{(2)}_{\Phi}$ and the first moment measure $\mu^{(1)}_{\Phi}$. Since $\alpha^{(2)}_{\Phi}$ can be computed with a second-order derivative of the PGFl similarly to \eqref{EqFirstMomentComputation}, $\mu^{(2)}_{\Phi}$ can be expressed as a combination of differentiated PGFls. While a relation between factorial and non factorial moments such as \eqref{EqMomentMeasureDecomp1c} exists for higher orders, it becomes increasingly complicated and tedious to write in order to produce the expression of non factorial moment measures.\newline
      \\
      If we compare the expression of factorial and non factorial moment measures in \eqref{EqMomentMeasureDecomp1c}, we can see that $\alpha^{(2)}_{\Phi}(B_1, B_2)$ is the joint expectation of \textit{distinct} points of the point process falling in $B_1$ and $B_2$. When $B_1$ and $B_2$ are disjoint, of course, factorial and non factorial moments are equivalent and $\alpha^{(2)}_{\Phi}(B_1, B_2)$ has the same interpretation as $\mu^{(2)}_{\Phi}(B_1, B_2)$. When $B_1$ and $B_2$ are \textit{not} disjoint, however, the factorial moment $\alpha^{(2)}_{\Phi}$ fails to take into account the joint occurrence of a single target in both regions $B_1, B_2$ and the quantity $\alpha^{(2)}_{\Phi}(B_1, B_2)$ is difficult to interpret.\newline
      \\
      So why can't we write $\mu^{(2)}_{\Phi}(B_1, B_2)$ directly as the derivative of a \textit{single} PGFl? Since we need to consider the occurrence of a single point of the point process falling in the intersection of the two regions, the quantity
      \begin{equation} \label{EqProductIndicator}
	1_{B_1}(x)1_{B_2}(x),
      \end{equation}
      among others, must appear somewhere during the derivation process. Let us differentiate twice the PGFl in the directions $1_{B_1}(.)$ and $1_{B_2}(.)$ and see what we get. Recall from previous lectures (see \eqref{eq:defGenFunc}) that the PGFl of a process $\Phi$ is given by
      \begin{equation} \label{EqPGFlDefinition}
	G_{\Phi}(h) = \sum_{n \geqslant 0} \int_{\Xcal^{(n)}} \left(\prod_{i = 1}^n h(x_i)\right) P_{\Phi}(d\{x_1,\dots,x_n\}).
      \end{equation}
      Thus, the differentiated PGFl reads
      \begin{multline} \label{EqPGFlDerivation1}
	\left.\delta^2 G_{\Phi}(h; 1_{B_1}, 1_{B_2})\right|_{h = 1}
	\\
	= \left.\delta^2 \left(\sum_{n \geqslant 0} \int_{\Xcal^{(n)}} \left(\prod_{i = 1}^n \cdot(x_i)\right) P_{\Phi}(d\{x_1,\dots,x_n\})\right)(h; 1_{B_1}, 1_{B_2})\right|_{\mathclap{~~~~h = 1}}.
      \end{multline}
      The expression \eqref{EqPGFlDerivation1} indicates that:
      \begin{enumerate}
	\item The functional to differentiate is the function on test functions $t$ defined by: $t \mapsto \ds{\sum_{n \geqslant 0} \int_{\Xcal^{(n)}}} \left(\prod_{i = 1}^n t(x_i)\right) P_{\Phi}(d\{x_1,\dots,x_n\})$
	\item The functional is differentiated in directions (or increments) $1_{B_1}$ and $1_{B_2}$;
	\item The differentiated functional is evaluated at the constant test function $h = 1$, i.e.: $\forall x \in \Xcal, ~h(x) = 1$.
      \end{enumerate}
      However, \eqref{EqPGFlDerivation1} is rather cumbersome and we will favour the slightly more compact notation
      \begin{multline} \label{EqPGFlDerivation2}
	\left.\delta^2 G_{\Phi}(h; 1_{B_1}, 1_{B_2})\right|_{h = 1} 
	\\
	= \left.\delta^2 \left(\sum_{n \geqslant 0} \int_{\Xcal^{(n)}} \left(\prod_{i = 1}^n h(x_i)\right) P_{\Phi}(d\{x_1,\dots,x_n\}); 1_{B_1}, 1_{B_2}\right)\right|_{\mathclap{~~~~h = 1}},
      \end{multline}
      where we must keep in mind that the argument of the PGFl is the test function $h$.\newline
      \\
      Since the sum and the integral are continuous linear operators, we can move the differentiation operator within the integral\footnotemark[2] and we get\footnotetext[2]{See exercise 1 in section \ref{SecExercise}}
      \begin{multline} \label{EqPGFlDerivation3}
	\left.\delta^2 G_{\Phi}(h; 1_{B_1}, 1_{B_2})\right|_{h = 1}
	\\
	= \sum_{n \geqslant 0} \int_{\Xcal^{(n)}} \left.\delta^2 \left(\prod_{i = 1}^n h(x_i); 1_{B_1}, 1_{B_2}\right)\right|_{h = 1} P_{\Phi}(d\{x_1,\dots,x_n\}).
      \end{multline}
      Let us expand a derivation term in \eqref{EqPGFlDerivation3}, e.g.\ for $n = 3$. If we differentiate $h(x_1) h(x_2) h(x_3)$ once in the direction $1_{B_1}$ we get:
      \begin{subequations} \label{EqPGFlDerivation4}
	\begin{align}
	  &\delta(h(x_1) h(x_2) h(x_3); 1_{B_1}) \nonumber
	  \\
	  &= \delta(h(x_1); 1_{B_1}) h(x_2) h(x_3) + h(x_1) \delta(h(x_2); 1_{B_1}) h(x_3) \nonumber
	  \\
	  &\quad\quad\quad\quad\quad\quad\quad\quad\quad\quad\quad\quad\quad\quad\quad\quad+ h(x_1) h(x_2) \delta(h(x_3); 1_{B_1}) \label{EqPGFlDerivation4a}
	  \\
	  &= 1_{B_1}(x_1) h(x_2) h(x_3) + h(x_1) 1_{B_1}(x_2) h(x_3) + h(x_1) h(x_2) 1_{B_1}(x_3) \label{EqPGFlDerivation4b}
	\end{align}
      \end{subequations}
      We see in \eqref{EqPGFlDerivation4b} that $h$ is not ``available'' at a given point $x_i$ once it has been differentiated -- for example, the first term does not contain $h(x_1)$ anymore. If we differentiate a second time in the direction $1_{B_2}$ then set $h = 1$ we get:
      \begin{align}
	&\left.\delta^2(h(x_1) h(x_2) h(x_3); 1_{B_1}, 1_{B_2})\right|_{h = 1} \nonumber
	\\
	&= 1_{B_1}(x_1) 1_{B_2}(x_2) + 1_{B_1}(x_1) 1_{B_2}(x_3) + 1_{B_2}(x_1) 1_{B_1}(x_2) \nonumber
	\\
	&\quad\quad\quad\quad\quad\quad+ 1_{B_1}(x_2) 1_{B_2}(x_3) + 1_{B_2}(x_1) 1_{B_1}(x_3) + 1_{B_2}(x_2) 1_{B_1}(x_3) \label{EqPGFlDerivation5}
      \end{align}
      Because the test function $h$ ``disappeared'' in the simple derivation terms such as $\delta(h(x_1); 1_{B_1})$, the desired products such as $1_{B_1}(x_1)1_{B_2}(x_1)$ do not appear in \eqref{EqPGFlDerivation5} -- in other words, the PGFl is not adapted to the production of \textit{non factorial} moment measures.\newline
      \\
      Let us consider the transformation $h \rightarrow e^{-f}$ and consider $f$ as the new test function. One can show that
      \begin{equation} \label{EqDerivationExponentialFunctional}
	\delta(e^{-f(x_1)}; 1_{B_1}) = \delta(-f(x_1); 1_{B_1}) e^{-f(x_1)} = - 1_{B_1}(x_1) e^{-f(x_1)}.
      \end{equation}
      The result above may look familiar; indeed, if we consider functions rather than functional, we have the well-known result
      \begin{equation} \label{EqDerivationExponentialFunction}
	\left(e^{-f(x_1)}\right)' = -f'(x) e^{-f(x)}.
      \end{equation}
      This time, the test function $f$ did \textit{not} disappear in the derivation process \eqref{EqDerivationExponentialFunctional}; differentiating a second time in direction $1_{B_2}$ yields
      \begin{equation} \label{EqDerivationExponentialFunctional2}
	\delta^2(e^{-f(x_1)}; 1_{B_1}, 1_{B_2}) = - 1_{B_1}(x_1) \delta(e^{-f(x_1)}, 1_{B_2}) = 1_{B_1}(x_1) 1_{B_2}(x_1) e^{-f(x_1)},
      \end{equation}
      and the desired product $1_{B_1}(x_1) 1_{B_2}(x_1)$ \textit{do} appear.\newline
      \\
      More generally, the \textit{Laplace functional} \eqref{eq:LaplaceFunctional} of the process
      \begin{equation} \label{EqLaplaceDefinition}
	L_{\Phi}(f) = G_{\Phi}(e^{-f}) = \sum_{n \geqslant 0} \int_{\Xcal^{(n)}} \left(\prod_{i = 1}^n e^{-f(x_i)}\right) P_{\Phi}(d\{x_1,\dots,x_n\})
      \end{equation}
      is well-adapted to the production of \textit{non factorial} moment measures since one can show that
      \begin{multline} \label{EqLaplaceMoment}
	\forall B_1,B_2,\dots,B_n \in \Bcal(\Xcal),
	\\
	\mu^{(n)}_{\Phi}(B_1,B_2,\dots,B_n) = (-1)^n \left.\delta^n L_{\Phi}(f; 1_{B_1},\dots,1_{B_n})\right|_{f = 0}.
      \end{multline}  
      In particular, the quantity $\mu^{(2)}_{\Phi}(B, B)$ in the expression of the variance \eqref{EqVariance} is given by the second-order derivative $\left.\delta^2 L_{\Phi}(f; 1_{B},1_{B})\right|_{f = 0}$.
      
    \section{Example: PHD filter with variance} \label{SecPHD}
      Let us have a look at what we can do with the statistics $(\mu^{(1)}_{\Phi}(\cdot), var_{\Phi}(\cdot))$.\newline
      \\
      Recall that the principle of the PHD filter is to propagate the \textit{first moment density} or \textit{intensity} or \textit{Probability Hypothesis Density} $D_{\Phi}(x)$, whose integral in any region $B \in \Bcal(\Xcal)$ yields the mean target number $\mu^{(1)}_{\Phi}(B)$. More formally, the PHD $D_{\Phi}$ is the Radon-Nikodym derivative of the first moment measure $\mu^{(1)}_{\Phi}$ w.r.t. the Lebesgue measure:
      \begin{equation} \label{EqPHD}
	\forall B \in \Bcal(\Xcal),~\mu^{(1)}_{\Phi}(B) = \int_B D_{\Phi}(x)dx
      \end{equation}
      The PHD filtering mechanisms can be described by considering moment measures or the corresponding moment densities. Since it is unclear whether the variance does possess such a Radon-Nikodym derivative w.r.t. the Lebesgue measure, we will always write the variance as a function on the $\sigma$-algebra $\Bcal(\Xcal)$; that is, we will always evaluate the variance of the process in a \textit{region} of the state space $\Xcal$, not at a single point. For the sake of consistency, the first moment will be expressed through its \textit{measure} $\mu^{(1)}_{\Phi}$ rather than its \textit{density} $D_{\Phi}$.\newline
      \\
      The PHD filtering mechanism can then be illustrated for one time step:
      \begin{figure}[H]
	\centering
	\includegraphics[width=\textwidth]{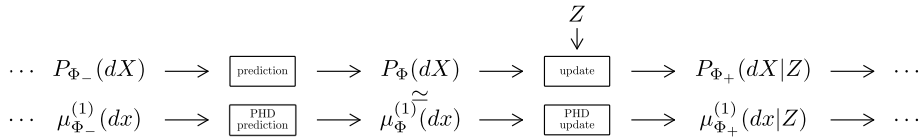}
	\caption{PHD filtering: data flow\label{FigFilterProcess}}
      \end{figure}
      The Poisson assumption states that the predicted process ($\Phi$ in figure \ref{FigFilterProcess}) is a Poisson process with a first moment measure equal to the output of the PHD prediction step $\mu^{(1)}_{\Phi}$. Because the prior process ($\Phi_{-}$ in figure \ref{FigFilterProcess}) is not necessarily Poisson, neither is the predicted process $\Phi$ after the full multi-object Bayes prediction. In other words, the probability distribution $P_{\Phi}$ is \textit{not necessarily} the probability distribution of a Poisson process, and the Poisson assumption is therefore an \textit{approximation} -- necessary to close the PHD recursion as exposed during the tutorials.\newline
      \\
      Now that second-order information on the process is available through its variance, we can ask ourselves the following questions:\newline
      \\
      \indent Q1: ``Let us assume that the predicted process $\Phi$ is Poisson. We know that scalar Poisson distributions are completely determined by their mean $\lambda$ -- in particular, the variance of a scalar Poisson distribution is also $\lambda$. Does this result extend to point processes, i.e.\ is the variance $var_{\Phi}(B)$ equal to the first moment measure $\mu^{(1)}_{\Phi}(B)$ for any region $B \in \Bcal(\Xcal)$ \textit{if $\Phi$ is a Poisson point process}?''\newline
      \\
      \indent Q2: ``What is the variance $var_{\Phi_{+}}(B|Z)$ of the updated process $\Phi_{+}$ in some region $B \in \Bcal(\Xcal)$? How is this variance related to the amount of information provided by the current set of measurements $Z$?''\newline
      \\
      The first question will now be addressed, while the second one is left in exercise \ref{SubsecPHDVariance}.\newline
      \\
      Let us assume that the predicted process $\Phi$ is Poisson, with first moment measure $\mu^{(1)}_{\Phi}$ (output of the PHD prediction step). Recall from the previous lectures, e.g.\ \eqref{eq:PoissonPgfl}, that the PGFl of a Poisson process is given by
      \begin{equation} \label{EqPoissonPGFl}
	G_{\Phi}(h) = \exp\left[\int (h(x) - 1)\mu^{(1)}_{\Phi}(dx)\right].
      \end{equation}
      Since our goal is to find the expression of the variance $var_{\Phi}$, following the definition \eqref{EqVariance} we need to find the expression of the second moment measure $\mu^{(2)}_{\Phi}$ first, and we have just seen (section \ref{SecLaplace}) that we can compute this quantity from the Laplace functional of the process. From the general expression of the Laplace functional \eqref{EqLaplaceDefinition} and the PGFl of a Poisson point process \eqref{EqPoissonPGFl} follows the Laplace functional of a Poisson point process
      \begin{equation} \label{EqLaplacePGFl}
	L_{\Phi}(f) = \exp\left[\int (e^{-f(x)} - 1)\mu^{(1)}_{\Phi}(dx)\right].
      \end{equation}
      Let us fix some region $B \in \Bcal(\Xcal)$. Using the relation between moment measures and derivatives of the Laplace functional \eqref{EqLaplaceMoment} we can write:
      \begin{subequations} \label{EqPoissonSecondMoment1}
	\begin{align}
	  \mu^{(2)}_{\Phi}(B,B) &= \left. \delta^2 L_{\Phi}(f; 1_B, 1_B) \right|_{f = 0} \label{EqPoissonSecondMoment1a}
	  \\
	  &= \left. \delta^2 \left(\exp\left[\int (e^{-\cdot(x)} - 1)\mu^{(1)}_{\Phi}(dx)\right]\right)(f; 1_B, 1_B) \right|_{f = 0} \label{EqPoissonSecondMoment1b}
	  \\
	  &= \left. \delta^2 \left(\exp \circ L\right)(f; 1_B, 1_B) \right|_{f = 0}, \label{EqPoissonSecondMoment1c}
	\end{align}
      \end{subequations}
      where the inner functional $L$ is defined as $L(f) = \int (e^{-f(x)} - 1)\mu^{(1)}_{\Phi}(dx)$. The derivation in \eqref{EqPoissonSecondMoment1c} can be easily expanded with Fa\`{a} di Bruno's formula for functional derivatives. Here the partitioning of the increments $1_B, 1_B$ is particularly simple: a single one-element partition $\{1_B, 1_B\}$, a single two-element partition $\{1_B\} - \{1_B\}$. Thus, \eqref{EqPoissonSecondMoment1c} becomes
      \begin{multline} \label{EqPoissonSecondMoment2}
	\mu^{(2)}_{\Phi}(B,B)
	\\
	= \left. \delta^2 \exp (L(f); \delta L(f; 1_B), \delta L(f; 1_B)) \right|_{f = 0} + \left.\delta \exp (L(f); \delta^2 L(f; 1_B, 1_B)) \right|_{f = 0}.
      \end{multline}
      Let us have a look at the increment $\delta L(f; 1_B)$ in \eqref{EqPoissonSecondMoment2}. Since the integral is a linear operator \footnotemark[2], we can move the derivation inside the integral, i.e.:\footnotetext[2]{See exercise 1 in section \ref{SecExercise}}
      \begin{subequations} \label{EqPoissonSecondMoment3}
	\begin{align}
	  \delta L(f; 1_B) &= \delta \left(\int (e^{-f(x)} - 1)\mu^{(1)}_{\Phi}(dx); 1_B\right) \label{EqPoissonSecondMoment3a}
	  \\
	  &= \int \delta(e^{-f(x)} - 1; 1_B)\mu^{(1)}_{\Phi}(dx) \label{EqPoissonSecondMoment3b}
	  \\
	  &= \int \delta(e^{-f(x)}; 1_B)\mu^{(1)}_{\Phi}(dx) \label{EqPoissonSecondMoment3c}
	\end{align}
      \end{subequations}
      Then, using the previously established rule \eqref{EqDerivationExponentialFunctional}:
      \begin{equation}
	\delta L(f; 1_B) = \int (-1_B(x)) e^{-f(x)} \mu^{(1)}_{\Phi}(dx). \label{EqPoissonSecondMoment4}
      \end{equation}
      As expected, the test function $f$ has not ``disappeared'' and thus $\delta L(f; 1_B)$ can be differentiated \textit{another time} without vanishing. Following the same reasoning that led to \eqref{EqPoissonSecondMoment4}, we can write the expression of the second-order derivative of the inner function $L$:
      \begin{subequations} \label{EqPoissonSecondMoment5}
	\begin{align}
	  \delta^2 L(f; 1_B, 1_B) &= \int (-1_B(x))^2 e^{-f(x)} \mu^{(1)}_{\Phi}(dx) \label{EqPoissonSecondMoment5a}
	  \\
	  &= \int 1_B(x) e^{-f(x)} \mu^{(1)}_{\Phi}(dx). \label{EqPoissonSecondMoment5b}
	\end{align}
      \end{subequations}
      Now that the increments $\delta L(f; 1_B)$ and $\delta^2 L(f; 1_B, 1_B)$ are known, let us resolve the outer differentiation in \eqref{EqPoissonSecondMoment2}. The outer functionals being exponentials, a nice rule similar to \eqref{EqDerivationExponentialFunctional} can be applied to resolve them; of course, going back to the definition of the chain rule would lead to the same result. The rule states that
      \begin{equation} \label{EqDerivationExponentialFunctionalExtension}
	\delta \exp (F(h); \delta F(h; \eta)) = \delta F(h; \eta) e^{F(h)}.
      \end{equation}
      Again, the functional derivation process of an exponential is similar to what we are used to with classical derivatives. Applying the new rule \eqref{EqDerivationExponentialFunctionalExtension} to \eqref{EqPoissonSecondMoment2} gives
      \begin{multline} \label{EqPoissonSecondMoment6}
	\mu^{(2)}_{\Phi}(B,B)
	\\
	= \left.\left(\delta L(f; 1_B)\right)^2 \exp (L(f))\right|_{f = 0} + \left.\delta^2 L(f; 1_B, 1_B)) \exp (L(f)) \right|_{f = 0}.
      \end{multline}
      We can then substitute in \eqref{EqPoissonSecondMoment6} the values of the increments which have been determined in \eqref{EqPoissonSecondMoment4} and \eqref{EqPoissonSecondMoment5b}:
      \begin{subequations} \label{EqPoissonSecondMoment7}
	\begin{align}
	  &\mu^{(2)}_{\Phi}(B,B) \nonumber
	  \\
	  &= \left.\left(\int (-1_B(x)) e^{-f(x)} \mu^{(1)}_{\Phi}(dx)\right)^2\right|_{f = 0} \left.\exp (L(f)) \right|_{f = 0} \nonumber
	  \\
	  &\quad\quad\quad\quad\quad\quad\quad\quad+ \left.\left(\int 1_B(x) e^{-f(x)} \mu^{(1)}_{\Phi}(dx)\right)\right|_{f = 0} \left.\exp (L(f)) \right|_{f = 0} \label{EqPoissonSecondMoment7a}
	  \\
	  &= \left(-\int_B \mu^{(1)}_{\Phi}(dx)\right)^2 \exp (L(0)) + \left(\int_B \mu^{(1)}_{\Phi}(dx)\right) \exp (L(0)) \label{EqPoissonSecondMoment7c}
	  \\
	  &= \left(-\mu^{(1)}_{\Phi}(B)\right)^2 \cdot 1 + \left(\mu^{(1)}_{\Phi}(B)\right) \cdot 1 \label{EqPoissonSecondMoment7d}
	  \\
	  &= \left(\mu^{(1)}_{\Phi}(B)\right)^2 + \mu^{(1)}_{\Phi}(B). \label{EqPoissonSecondMoment7e}
	\end{align}
      \end{subequations}
      Now that the second moment measure has been determined in \eqref{EqPoissonSecondMoment7e}, we can produce the variance $var_{\Phi}$ of the Poisson process $\Phi$ from the definition \eqref{EqVariance}:
      \begin{subequations} \label{EqVariancePoisson1}
	\begin{align}
	  var_{\Phi}(B) &= \mu^{(2)}_{\Phi}(B, B) - \left(\mu^{(1)}_{\Phi}(B)\right)^2 \label{EqVariancePoisson1a}
	  \\
	  &= \left(\mu^{(1)}_{\Phi}(B)\right)^2 + \mu^{(1)}_{\Phi}(B) - \left(\mu^{(1)}_{\Phi}(B)\right)^2 \label{EqVariancePoisson1b}
	  \\
	   &= \mu^{(1)}_{\Phi}(B). \label{EqVariancePoisson1c}
	\end{align}
      \end{subequations}
      As it may be expected, the variance \textit{of a Poisson process} equals its mean in any region $B \in \Bcal(\Xcal)$.\newline
      \\
      An interesting consequence to \eqref{EqVariancePoisson1c} is that the variance \textit{of a Poisson process} in any region $B$ is bounded by the mean target number in the state space since
      \begin{equation} \label{EqVarianceBoundsPoisson}
	var_{\Phi}(B) = \mu^{(1)}_{\Phi}(B) \leqslant \mu^{(1)}_{\Phi}(\Xcal).
      \end{equation}
      This means that the local behaviour of a Poisson process with ``reasonable global average behaviour'' -- i.e.\ with a finite mean target number in the whole state space $\mu^{(1)}_{\Phi}(\Xcal)$ -- can be estimated in any region $B$ with ``some accuracy'' since the variance of the target number in $B$ is finite as well according to \eqref{EqVarianceBoundsPoisson}. We will see in an exercise \ref{SubsecIIDProcess} that this is not necessarily true for more advanced (i.e.\ non Poisson) point processes.
      
    \section{Exercises} \label{SecExercise}
      \subsection{Functional derivatives of linear functionals} \label{SubsecLinearFunctional}
	A linear functional is a functional $L$ such that for any functions $f$, $g$, and for any scalars $\alpha$, $\beta$:
	\begin{equation} \label{ExLinearFunctional}
	  L(\alpha f + \beta g) = \alpha L(f) + \beta L(g).
	\end{equation}
	Let $L$ be a continuous linear functional. Prove that, for any functional $F$, any test function $h$ and any admissible direction $\eta$:
	\begin{equation} \label{ExLinearDerivative}
	  \delta (L \circ F)(h; \eta) = L(\delta F(h; \eta)).
	\end{equation}     
	\begin{hint}
	  Use the definition of the chain differential.
	\end{hint}
	This result allows us to ``switch'' integral and differentiation operators. Indeed, let us choose $L$ as
	\begin{equation} \label{ExLinearExample}
	  L(F(h)) = \int F(h)(x)\mu(dx)
	\end{equation}
	where $\mu$ is some measure on the $\sigma$-algebra $\Bcal(\Xcal)$ - for example, the first moment measure of some point process $\Phi$. Using \eqref{ExLinearDerivative} we can write:
	\begin{subequations} \label{ExLinearIntegral1}
	  \begin{align} 
	    \delta \left(\int F(\cdot)(x)\mu(dx)\right)(h; \eta) &= \delta (L \circ F)(h; \eta) \label{ExLinearIntegral1a}
	    \\
	    &= L(\delta F(h; \eta)) \label{ExLinearIntegral1b}
	    \\
	    &= \int \delta \left(F(\cdot)\right)(h(x); \eta) \mu(dx). \label{ExLinearIntegral1c}
	  \end{align} 
	\end{subequations}
	
      \subsection{PGFl and Laplace functional} \label{SubsecLinearPGFlLaplace}
	Let $\Phi$ be a point process with PGFl $G_{\Phi}$ and Laplace functional $L_{\Phi}$. Prove that:
	\begin{equation} \label{ExPGFlandLaplace}
	  \forall B \in \Bcal(\Xcal),~\left.\delta G_{\Phi}(h; 1_B)\right|_{h = 1} = -\left.\delta L_{\Phi}(f; 1_B)\right|_{f = 0}.
	\end{equation}
	What can we conclude on the first moment $\mu^{(1)}_{\Phi}$ and first factorial moment $\alpha^{(1)}_{\Phi}$?
	
      \subsection{PHD update and variance} \label{SubsecPHDVariance}
	Show that the variance of the updated process ($\Phi_{+}$ in figure \ref{FigFilterProcess}) is
	\begin{multline} \label{ExPHDUpdate}
	  var_{\Phi_{+}}(B|Z) = \int_B L(\phi|x)\mu^{(1)}_{\Phi}(dx)
	  \\
	  + \sum_{z \in Z}\frac{\int_B L(z|x)\mu^{(1)}_{\Phi}(dx)}{c(z) + \int L(z|x)\mu^{(1)}_{\Phi}(dx)}\left(1 - \frac{\int_B L(z|x)\mu^{(1)}_{\Phi}(dx)}{c(z) + \int L(z|x)\mu^{(1)}_{\Phi}(dx)}\right),
	\end{multline}
	where:
	\begin{itemize}
	  \item $Z = \{z_1,\dots,z_m\}$ is the set of current observations;
	  \item $L(\phi|x) = 1 - p_d(x)$, where $p_d$ is the probability of detection;
	  \item $L(z|x) = p_d(x)\hat{L}(z|x)$, where $\hat{L}$ is the single-measurement / single-target likelihood function;
	  \item $c(z)$ is the intensity of the false alarm process, assumed Poisson.
	\end{itemize}
	
	\begin{hint}
	  We have seen in tutorial that the PGFl of the updated process is given by
	  \begin{multline} \label{ExPHDUpdatePGFl}
	    G_{\Phi_{+}}(h|Z)
	    \\
	    = \exp \left[\int (h(x) - 1)L(\phi|x)\mu^{(1)}_{\Phi}(dx)\right] \prod_{z \in Z} \frac{c(z) + \int h(x)L(z|x)\mu^{(1)}_{\Phi}(dx)}{c(z) + \int L(z|x)\mu^{(1)}_{\Phi}(dx)}.
	  \end{multline}
	  From \eqref{ExPHDUpdatePGFl}, determine the expression of the Laplace functional $L_{\Phi_{+}}$ and differentiate it twice in direction $1_B$ to get the expression of $var_{\Phi_{+}}(B|Z)$.
	\end{hint}
	
      \subsection{Variance of i.i.d. processes}  \label{SubsecIIDProcess}
	An independent and identically distributed (i.i.d.) process $\Phi$ is an extension of a Poisson process which plays an important role in the construction of the CPHD filter. It is characterized by:
	\begin{itemize}
	  \item Its first moment measure $\mu^{(1)}_{\Phi}$ or, equivalently, its first moment density $D_{\Phi}$: each target is i.i.d. according to the normalised intensity $\frac{D_{\Phi}(.)}{\int D_{\Phi}(x)dx}$;
	  \item Its cardinality distribution $\rho_{\Phi}$ such that $\sum_{n \geqslant 1} n \rho_{\Phi}(n) = \mu^{(1)}_{\Phi}(\Xcal)$: $\rho_{\Phi}(n)$ is the probability that there are exactly $n$ targets in the state space $\Xcal$.
	\end{itemize}
	The PGFl of a i.i.d. point process is
	\begin{equation} \label{ExIIDPGFl}
	  G_{\Phi}(h) = \sum_{n \geqslant 0} \rho_{\Phi}(n) \left(\frac{\int h(x) \mu^{(1)}_{\Phi}(dx)}{\mu^{(1)}_{\Phi}(\Xcal)}\right)^n.
	\end{equation}
	Let $\Phi$ be a i.i.d. process, with first moment measure $\mu^{(1)}_{\Phi}$ and cardinality distribution $\rho_{\Phi}$.\newline
	\\
	a) Prove that in the special case where the cardinality distribution is Poisson with parameter $\mu^{(1)}_{\Phi}(\Xcal)$, $\Phi$ is a Poisson process.
	
	\begin{hint}
	  Show that if $\rho_{\Phi}(n) = e^{-\mu^{(1)}_{\Phi}(\Xcal)} \frac{\left(\mu^{(1)}_{\Phi}(\Xcal)\right)^n}{n!}$, then the PGFl \eqref{ExIIDPGFl} reduces to the PGFl of a Poisson process \eqref{EqPoissonPGFl}.
	\end{hint}
	
	b) Show that the variance of $\Phi$ in any region $B \in \Bcal(\Xcal)$ is equal to
	\begin{equation}
	  var_{\Phi}(B) = \mu^{(1)}_{\Phi}(B) + \left(\mu^{(1)}_{\Phi}(B)\right)^2\left[\frac{\sum_{n \geqslant 2} n(n - 1) \rho_{\Phi}(n) }{\left[\sum_{n \geqslant 1} n \rho_{\Phi}(n) \right]^2} - 1\right].
	\end{equation}
	What happens if the cardinality distribution is Poisson?\newline
	\\
	c) \textit{(Advanced)} i.i.d. processes have a less intuitive behaviour than Poisson processes and can yield surprising results.\newline
	\\
	Prove that, for any mean target number $0 < \mu < \infty$, any constant $0 < C < \infty$, there exists a i.i.d. point process $\Phi_C$ such that:
	\begin{itemize}
	  \item $\int \mu^{(1)}_{\Phi_C}(dx) = \mu$;
	  \item $\forall B \in \Bcal(\Xcal),~var_{\Phi_C}(B) \geqslant \mu^{(1)}_{\Phi_C}(B) + C \cdot \lambda(B)$;
	\end{itemize}
	where $\lambda$ is the Lebesgue measure.
	
	\begin{hint}
	  Produce a sequence of i.i.d. point processes $\{\Phi_s\}_{s \geqslant 0}$ such that:
	  \begin{itemize}
	  \item $\forall s \in \Nset,~\int \mu^{(1)}_{\Phi_s}(dx) = \int \mu^{(1)}_{\Phi}(dx)\mu$;
	  \item $\forall B \in \Bcal(\Xcal),~\lim_{s \rightarrow \infty} \frac{var_{\Phi_s}(B) - \mu^{(1)}_{\Phi_s}(B)}{\lambda(B)} = \infty$.
	  \end{itemize}
	\end{hint}
	That is, contrary to a Poisson process (see section \ref{SecPHD}), one can always find a i.i.d. process $\Phi$ with a ``reasonable average behaviour'' -- i.e.\ a finite \textit{global} mean target number $\mu^{(1)}_{\Phi}(\Xcal)$ -- and yet an \textit{arbitrary high} variance in a given region $B$ -- the information on the \textit{local} target number in $B$ is ``completely unreliable''.
  
\bibliography{NotesSummerSchool.bib}
\bibliographystyle{plain}

\appendix
\chapter[Review material]{Review material}

\section{Probability theory, random variables and differentiation}

The following supplementary review material is taken from the
Erasmus Mundus Vision and Robotics course 
taught by Daniel Clark at Heriot-Watt, from 2012, 
and at the First International Summer School on Finite Set Statistics 
in July 2013.

\subsection{Set theory}

\begin{definition}[Sample space]
The set $S$ of all possible outcomes of
a particular experiment is called the sample space.
\end{definition}

\begin{example}[Coin tossing]
In the experiment of tossing a coin,
the sample space contains two outcomes:
\begin{align*}
S=\{H,T\}.
\end{align*}
\end{example}

\begin{example}[Waiting time]
Consider an experiment where the 
observation is the time that it takes 
for an electronic component to fail.
Then the sample space is all positive numbers,
\begin{align*}
S=(0,\infty).
\end{align*}
\end{example}

Sample spaces can be classified into two types,
according to the number of elements they contain.
They can either be {\it countable}
or {\it uncountable}.

\begin{definition}[Event]
An event is any collection of possible
outcomes of an experiment, that is,
any subset of sample space $S$.
Let $A$ be an event, so that $A\in S$.
We say that $A$ occurs is the outcome of the 
experiment is in the set $A$.
\end{definition}

\begin{definition}[Probability Axioms]

A probability measure satisfies the following:
\begin{align*}
p(A)&\ge 0,
\\
p(S)&=1,
\\
p(\emptyset)&=0,
\end{align*}
and, for disjoint $A_i$, $i>0$, we have
\begin{align*}
p(\cup_{i=1}^\infty A_i) = \sum_{i=1}^\infty p(A_i).
\end{align*}

\end{definition}

\subsection{Discrete random variables}

\begin{definition}[Random variable]
A random variable is a mapping from the sample
space into the real numbers, or vectors.
\end{definition}

For any discrete random variable $X$, the {\it mean} value
$\mathbb{E}[X]$ is an indication of the centre of the 
distribution of $X$.
This is also known as the first-order moment of $X$.

\begin{definition}[Expectation of a discrete random variable]
The expected value, mean, or expectation of a 
discrete random variable $X$,
written $\mathbb{E}[X]$
is given by 
\begin{align*}
\mathbb{E}[X]=
\sum_{x}
xP(X=x)
\end{align*}
\end{definition}

\begin{myrule}[Properties of expectation for random variables]
\begin{align*}
\mathbb{E}[c]&=c& \forall c\in\R\\
\mathbb{E}[cX]&=c\mathbb{E}[X] &\forall c\in\R\\
\mathbb{E}[X+Y]&=\mathbb{E}[X]+\mathbb{E}[Y]
\end{align*}
\end{myrule}

The $k^{th}$-order moment of $X$ is defined to be 
$\mathbb{E}[X^k]$ for $k=1,2,3,\ldots$.

\begin{definition}[Variance]
The 
variance of 
random variable $X$
is 
can be determined from the moments with
\begin{align*}
\var(X)=\mathbb{E}[(X-\mathbb{E}[X])^2]
\end{align*}
\end{definition}

\begin{myexercise}
Using the properties of expectation,
show that 
\begin{align*}
\var(X)=\mathbb{E}[X^2]-\mathbb{E}[X]^2.
\end{align*}
\end{myexercise}

\newpage

\subsection{Ordinary  derivatives}

In this section we revise some 
of the rules for ordinary differentiation
for product of functions 
and composite functions.
We shall extend these rules to 
the calculus of variations in the following section.
We begin with the basic notions of ordinary and partial
derivative.
Ordinary derivatives can be viewed as
a special case of partial derivative 
for one variable.

\begin{definition}[Ordinary derivative]
The function $f:X\rightarrow Y$, where $X=Y=\mathbb{R}$, has a 
derivative if the following limit exists
\begin{equation*}
{d \over  dx} f(x) = \lim_{\epsilon\rightarrow 0}
\dfrac{1}{\epsilon} \left( f(x+\epsilon)-f(x) \right).
\end{equation*}
\end{definition}

\begin{myrule}[Product rule]
The product rule for ordinary derivatives is
\begin{align*}
{d\over dx}
\left(
f(x)g(x)
\right)
=f'(x)g(x)+f(x)g'(x)
\end{align*}
\end{myrule}

\begin{myrule}[Chain rule]
The chain rule for ordinary derivatives is
\begin{align*}
{d\over dx}
f\circ g(x)
=
f'(g(x))g'(x)
\end{align*}

\end{myrule}

\begin{myrule}[Leibniz's formula]
The higher-order product rule, known as
{\it Leibniz's formula},
is given by
\begin{align*}
{d^n
\over 
dx^n}
f(x)g(x)=
\sum_{k=0}^n
{n!\over (n-k)!k!}
f^{(k)}(x)g^{(n-k)}(x)
\end{align*}
\end{myrule}

\begin{myrule}[Fa\`a di Bruno's formula]
The higher-order chain rule, known as 
{\it Fa\`a di Bruno's formula},
is 
\begin{align*}
{d^r
\over 
dx^r}
f(g(x))=
\sum_{\pi\in \Pi(1,\ldots,r)}
f^{(|\pi|)}(g(x))
\prod_{\gamma\in\pi}
g^{(|\gamma|)}(x),
\end{align*}
where $\Pi(1,\ldots,r)$
is the set of all partitions of the set $\{1,\ldots,r\}$,
and $|\gamma|$ is the size of cell $\gamma$.
\end{myrule}

\section{Probability generating functions}

\subsection{Generating functions}
A useful way of describing an infinite sequence $u_0,u_1,u_2,\ldots$
of real numbers is to write down the {generating function}
of the sequence, defined as
\begin{align*}
G(s)=\sum_{k=0}^\infty u_n s^n.
\end{align*}

This is an important concept that can be used to model
systems of multiple objects where there is uncertainty in the number of objects.
This lecture is mainly concerned with applications of this concept,
and how it can be used.

\begin{example}[Taylor (or Maclaurin) series]
The Taylor series of 
a real-valued function $f(x)$
that is infinitely differentiable in the neighbourhood of 
$0$ is the power series
\begin{align*}
f(x)\approx 
\sum_{n=0}^\infty
{f^{(n)}(0)
\over n!}
x^n,
\end{align*}
where $f^{(n)}(0)$ is the $n^{th}$ derivative of $f$ evaluated at $0$.
\end{example}

\begin{definition}[Probability generating function (p.g.f.)]
Many random variables take values in the set of non-negative integers.
Let $p_X(k)=P(X=k)$, for $k=0,1,2,\ldots$
be the mass function that satistfies
$p_X(k)\ge 0$ for all non-negative integers $k$,
and $\sum_{k=0}^\infty p_X(k)=1$.
Then the {\it probability generating function} (p.g.f.) of $X$ is the function $G_X(s)$ 
defined by
\begin{align*}
G_X(s)=\sum_{k=0}^\infty p_X(k) s^k.
\end{align*}
\end{definition}
From this definition, we can immediately see that $P(X=0)=G_X(0)$, and $G_X(1)=1$.
In the following, we shall identify some other properties 
of the p.g.f..
Firstly, we illustrate the concept with two different parameterisations 
that you may be familiar with, namely
the {\it Bernoulli distribution}
and the {\it Poisson distribution}.

\begin{example}[Bernoulli distribution]
Let us suppose that $P(X=k)=0$ for all $k>1$.
Thus there is either zero or one object.
Then the p.g.f. becomes
\begin{align*}
G_X(s)=p_X(0)+p_X(1)s.
\end{align*}
It follows that $G_X(1)=p_X(0)+p_X(1)=1$, and that $G_X(0)=p_X(0)$.
\end{example}

\begin{example}[Poisson distribution]
The {Poisson distribution}
is a distribution that can be uniquely characterised by its {\it rate} $\lambda>0$.
The rate gives both the mean and variance of the number of objects.
The p.m.f. in this case becomes
\begin{align*}
p_X(k)=P(X=k)={1\over k!}\lambda^k\exp(-\lambda),
\end{align*}
for $k=0,1,2,\ldots$.
The generating function is then
\begin{align*}
G_X(s)=\sum_{k=0}^\infty{1\over  k!}\lambda^k\exp(-\lambda)s^k.
\end{align*}
\end{example}

The Poisson distribution is attractive to use since
its entire p.m.f. can be specified by $\lambda$,
and it models well naturally occurring phenomena, such radioactive decay.

\begin{myexercise}[Poisson p.g.f.]
Use the Taylor expansion of the exponential function,
\begin{align*}
\exp(x)=
\sum_{k=0}^\infty {x^k\over k!},
\end{align*}
to show that the Poisson p.g.f.
simplifies to 
\begin{align*}
G_X(s)= \exp(\lambda(s-1)).
\end{align*}
\end{myexercise}

\subsection{Factorial moments}

The mean and variance are the first- and second-order {\it central moments}.
Another useful kind of moment is the factorial moment, defined as follows.
\begin{definition}[Factorial moment]
The $r^{th}$-order factorial moment is defined with
\begin{align*}
\mathbb{E}\left[ X(X-1)(X-2)\ldots (X-r+1) \right].
\end{align*}
\end{definition}

The factorial moments are particularly convenient to use 
as statistical descriptors of the probability distribution
that we wish to model, since they can be determined from the p.g.f. by differentiation,
which is shown in the following lemma.

\begin{mylemma}
The factorial moments can be determined from the p.g.f. as follows
\begin{align*}
G_X^{(r)}(1)=\mathbb{E}\left[ X(X-1)(X-2)\ldots (X-r+1) \right].
\end{align*}
\end{mylemma}


\begin{myexercise}[Poisson distribution]
Find the $r^{th}$-order factorial moment
of the Poisson distribution.
(Differentiate $r$-times and set $s=1$.)
\end{myexercise}

\begin{myexercise}
Show that it is possible
to recover the p.m.f. $p_X(k)=P(X=k)$
via differentiation.
Use this to determine the p.m.f. of the Poisson distribution.
\end{myexercise}

\subsection{Sums of independent random variables }

\begin{mytheorem}
If $X$ and $Y$ are independent non-negative integer-valued
random variables with generating functions 
$G_X(s)$ and $G_Y(s)$ respectively.
Then the generating function of the 
summation $X+Y$ is given by 
\begin{align*}
G_{X+Y}(s)=G_X(s)G_Y(s).
\end{align*}
\end{mytheorem}

\begin{myproof}
\begin{align*}
G_{X+Y}(s)&=\mathbb{E}[s^{X+Y}]
\\
&=\mathbb{E}[s^Xs^Y]
\\
&=\mathbb{E}[s^X]\mathbb{E}[s^Y]
\\
&=G_X(s)G_Y(s).
\end{align*}
The third line above follows from the fact that for independent random varables $X$ and $Y$,
we have
\begin{align*}
\mathbb{E}[f(X)g(Y)]
&=\mathbb{E}[f(X)]\mathbb{E}[g(Y)],
\end{align*}
for any functions $f$ and $g$.
\end{myproof}

\begin{myexercise}
Determine the factorial moments and p.m.f. of the 
product 
\begin{align*}
G_X(s)G_Y(s).
\end{align*}
\noindent
(Hint: Use Leibniz' formula.)

\end{myexercise}

\subsection{Branching processes}

In this section we discuss branching processes.
These will form the basis for our model
for multi-object dynamics for multi-target tracking algorithms
in the following section on probability generating functionals.

Let us consider a sequence $Z_0, Z_1, Z_2,\ldots$
of non-negative integer valued random variables.
$Z_n$ is interpreted as the number of objects in the $n^{th}$
generation of a population.
For convenience, it is often assumed that the $Z_0=1$.
The probability that an object 
existing in the $n^{th}$-generation
has $k$ children
in the ${n+1}^{th}$-generation is given by $p_k$.
Note that we assume that this probability is independent
of the generation number.
We further assume that each object 
reproduces independently of other objects.

More explicitly, let us consider the probability
generating function
\begin{align*}
G(s)=\sum_{k=0}^\infty
p_k s^k,
\end{align*}
where $|s|\le 1$, and $s$ is a complex variable.
Now consider the sequence of generating functions formed by composition, i.e.
\begin{align*}
G_0(s)&=s
\\
G_1(s)&=G(s)
\\
\ldots
\\
G_{n+1}(s)&=G_n(G(s)),
\end{align*}
for $n=0,1,2,\ldots$.

\begin{myexercise}
Suppose that the generating function at generation $n$
is given by $G_n(s)$,
and that of generating $n+1$ is given by
$G_{n+1}(s)=G_\gamma(s)G_n(G(s))$.
Compute the 
first-order factorial moment of $G_{n+1}(s)$.

\end{myexercise}

\subsection{Joint generating functions}
To use generating functions for inference problems,
we need to consider generating functions of two variables.
We will firstly define joint generating functions,
and then use them to determine the conditional p.m.f.
and factorial moments via differentiation.

\begin{definition}[Joint generating function]
The joint generating function of 
non-negative integer-values
random variables $Z$ and $X$
is defined with
\begin{align*}
G_{Z,X}(s,t)
=\mathbb{E}[t^{Z}s^{X}]
=\sum_{m,n=0}^\infty
p_{Z,X}(m,n)t^ms^n.
\end{align*}
\end{definition}

\begin{example}
The conditional p.m.f. $p_{X|Z}$ is given by:
\begin{align*}
p_{X|Z}(n|Z = m) = {p_{Z,X}(m,n) \over p_{Z}(m)} = {{1 \over n!} \left.{\partial^{m+n} \over \partial t^m\partial s^n} G_{Z,X}(t,s) \right|_{t=0, s=0} \over \left.{\partial^{m}\over \partial t^m} G_{Z,X}(t,1) \right|_{t=0}} 
\end{align*}
$p_{X|Z}(n|Z = m)$ is the probability that $X = n$ given that $Z = m$.\newline
\\
Similarly, the conditional nth-order factorial moment $\mu^{(n)}_{X|Z}$ is given by:
\begin{align*}
\mu^{(n)}_{X|Z}(Z = m) = {\left.{\partial^{m+n} \over \partial t^m\partial s^n} G_{Z,X}(t,s) \right|_{t=0, s=1} \over \left.{\partial^{m}\over \partial t^m} G_{Z,X}(t,1) \right|_{t=0}} 
\end{align*}
$\mu^{(n)}_{X|Z}(Z = m)$ is the nth-order factorial moment of $X$ given that $Z = m$. In particular, $\mu^{(1)}_{X|Z}(Z = m)$ is the mean value of $X$ given that $Z = m$.
\end{example}

\begin{myexercise}
Consider the joint generating function
\begin{align*}
G_{Z,X}(t,s)=G_\kappa(t)G_p(s G_L(t)),
\end{align*}
where
\begin{align*}
G_\kappa(t)&=\exp(\lambda(t-1))
\\
G_p(s)&=\exp(\mu(s-1))
\\
G_L(t)&=1-p_d+p_d\cdot t.
\end{align*}
Find the 
conditional first factorial moment $\mu^{(1)}_{X|Z}(Z = m)$.
\end{myexercise}

\begin{myexercise}
Suppose that the joint generating function
factorises as the product of two generating functions, as follows
\begin{align*}
G_{Z,X}(t,s)=G_{X}(s)G_{Z}(t).
\end{align*}
Show that the conditional p.m.f. $p_{X|Z}(n|m)$
is equal to $p_X(n)$.
What does this say about the dependencies between $X$ and $Z$?
\end{myexercise}

\newpage
\section{Probability generating functionals}

In this section we extend the ideas 
from the previous section
to deal with multi-object probability densities.
In order to do this,
we need the notion of a
functional, or variational, derivative.
We introduce directional 
derivatives next, and
then show
how these are used to 
determine joint probability densities and
factorial moment densities,
from probability generating
functionals, 
a generalisation of 
probability generating functions.
Following this, we
apply these concepts to
spatial branching processes
and for Bayesian estimation.

\subsection{Functional derivative}

Functional derivatives are used in physics for calculating
quantities related to many-particle systems
in quantum field theory.
In multi-target tracking, we use them to determine
practical algorithms.

\begin{definition}[Functional derivative]
Suppose that we have a functional 
$f(h)$, whose argument is a function of a single-object state variable
$h(x)$.
Then the functional derivative of $f(h)$
in direction $\varphi(x)$ is given by
\begin{align*}
\delta f(h; \varphi)
=
{\lim_{\epsilon\rightarrow 0}}
{f(h+\epsilon\varphi)-f(h)
\over
\epsilon
}.
\end{align*}
\end{definition}
Usually, we differentiate with respect to
some arbitrary function $\varphi(x)$,
and then set $\varphi(x)=\delta_y(x)$.

\begin{example}[Linear functional]
Consider a linear functional defined with
\begin{align*}
f(h)=\int h(x)s(x) dx.
\end{align*}
Then the functional derivative is found with
\begin{align*}
\delta f(h;\varphi)
&=
{\lim_{\epsilon\rightarrow 0}}
{\int (h(x)+\epsilon\varphi(x))s(x)dx-\int h(x)s(x)dx
\over
\epsilon
}
\\\notag
&=
{\lim_{\epsilon\rightarrow 0}}
{\epsilon \int \varphi(x)s(x)dx
\over
\epsilon}
\\\notag
&=
\int\varphi(x)s(x)dx.
\end{align*}
Now set $\varphi(x)=\delta_y(x)$,
then
\begin{align*}
\delta f(h;\delta_y)
&=
s(y).
\end{align*}
\end{example}
The functional derivative obeys the usual linearity property,
and product rule.
Note that the chain rule below has a bit of a strange form
since it does not factorise in the usual way.
\begin{align*}
\delta (f_1+f_2)(h;\varphi)
&=
\delta f_1(h;\varphi)
+\delta f_2(h;\varphi)
\\
\delta (f_1 \cdot f_2) (h;\varphi)
&=
\delta f_1(h;\varphi)f_2(h)+
f_1(h) \delta f_2(h;\varphi)
\\
\delta (f_1 \circ f_2) (h;\varphi)
&=
\delta f_1(h;\delta f_2(h;\varphi)).
\end{align*}

\begin{example}[Exponential of linear functional]
Let us consider an example of the chain rule
on the composition of the exponential function and
 a linear functional
\begin{align*}
\delta \exp(f(h);\varphi) & = \lim_{\epsilon\rightarrow 0} \epsilon^{-1} \left(\exp(f(h)+\epsilon f(\varphi))-\exp(f(h))\right) \\
& = \lim_{\epsilon\rightarrow 0} \epsilon^{-1} \exp(f(h)) \exp(\epsilon f(\varphi)-1) \\
& = f(\varphi)\exp(f(h)),
\end{align*}
where we use the fact that
\begin{align*}
\lim_{\epsilon\rightarrow 0}
\epsilon^{-1}
(\exp(\epsilon y)-1)
=y.
\end{align*}

\end{example}

The next theorem, derived recently, generalises Fa\`a di Bruno's formula
for general differentials.
\begin{mytheorem}[Fa\`a di Bruno's formula \cite{Clark2013}]
The generalised form for the $n^{th}$-order derivative is as follows.
\label{thm:chainRule}
    \begin{equation*}
\delta^n (f\circ g)(x; \eta_1,\ldots,\eta_n ) = \sum_{\pi\in \Pi(\eta_1,\ldots,\eta_n)}
\delta^{|\pi|}f \left( g(x); \delta^{|\omega|} g \left( x; \xi : \xi \in \omega \right):\omega\in \pi \right)
\end{equation*}
\end{mytheorem}

\subsection{Generating functionals}

The generating functional generalises the generating function by allowing for a formal power series in a functional variable $h$.

\begin{definition}[Generating functional]
Let $h:\mathbf{X} \rightarrow \mathbb{R}$ be a test function. Every functional $G$ continuous in the field of continuous functions can be represented by the expression
\begin{equation*}
G(h) = \sum_{n=0}^{\infty}{1\over n!} \int \left( \prod_{i=1}^{n} dx_i\, h(x_i) \right) g_n(x_1,\ldots,x_n),
\end{equation*}
where the functions $g_n:\mathcal{X}\rightarrow\mathbb{R}$ are a 
sequence of continuous functions associated to $G$ and independent of $h$. 
By convention, the first term $g_0(\emptyset)$, which is a constant, is included in the summation.
\end{definition}

\begin{example}[Correlation functions]
In quantum field theory,
physicists describe many-particle systems through the
composition of an exponential function
and a generating functional
of connected Feynman diagrams,
so that
\begin{align*}
Z(h)=\exp(W(h)),
\end{align*}
where
\begin{align*}
W(h)=\sum_{n=0}^\infty
{1\over n!}
\left(
\int \prod_{i=1}^n
h(x_i)dx_i
\right)
w_n(x_1,\ldots,x_n).
\end{align*}
In this case,
the correlation 
functions $z_n(x_1,\ldots,x_n)$
that describe $Z(h)$ 
are found with
\begin{align*}
\left.
\delta^n
Z(h; \delta_{x_1},\ldots,\delta_{x_n})
\right|_{h=0},
\end{align*}
where
we find that
\begin{align*}
z_1(x)&=w_1(x)
\notag
\\
z_2(x_1,x_2)&=w_2(x_1,x_2)+w_1(x_1)w_1(x_2)
\notag
\\
z_3(x_1,x_2,x_3)&=w_3(x_1,x_2,x_3) +w_1(x_1)w_2(x_2,x_3)+ w_1(x_2)w_2(x_1,x_3)\\
& + w_1(x_3)w_2(x_2,x_1)+ w_1(x_1)w_1(x_2)w_1(x_3)
\notag
\end{align*}

\end{example}

\begin{definition}[Probability generating functional]
The {\it probability generating functional}, or p.g.fl.,
of a probability distribution
$P$
is defined 
as the expectation value of
the symmetric function $\omega(x^n)=\prod_{i=1}^n h(x_i)$, so that
\begin{equation*}
G(h) = \sum_{n=0}^{\infty}{1\over n!} \int
\left( \prod_{i=1}^{n} dx_i\, h(x_i) \right) p_n(x_1,\ldots,x_n),
\end{equation*}
\end{definition}

The p.g.fl. 
is convenient for 
stochastic modelling with point processes.
For example, suppose that we have two independent point processes
specified by p.g.fl.s $G_1(\xi)$ and $G_2(\xi)$.
Then it is straightforward to show that the superposition 
of these processes is found with
\begin{align*}
G(\xi)=G_{1}(\xi) G_{2}(\xi).\label{superposition}
\end{align*}

\begin{example}[Bernoulli point process]
Let us suppose that $P(X=k)=0$ for all $k>1$.
Thus there is either zero or one object.
Then the p.g.f. becomes
\begin{align*}
G(h)=p_0+\int p_1(x) h(x)dx.
\end{align*}
\end{example}

\begin{example}[independent, identically distributed (i.i.d.) cluster  process]

\begin{equation*}
\label{eq:pgfl1}
G(h) = \sum_{n=0}^{\infty}{1\over n!} \int
\rho(n)\left( \int dx\, h(x) s(x) \right)^n,
\end{equation*}
where $\rho(n)$ is a probability mass function,
and each object is 
independent, and identically distributed 
according to $s(x)$.

\end{example}

\begin{example}[Poisson process]

The p.g.fl. for the Poisson process is given by
\begin{align*}
G(h)= \exp(\lambda(f(h)-1)),
\end{align*}
where $f(h)=\int h(x)s(x)dx$ is a linear functional.
\end{example}

\subsubsection*{Using the p.g.fl.}

In a similar way to finding the factorial moments
from the p.g.f., we can find factorial moment
densities from the p.g.fl. via differentiation.
However, instead of using ordinary derivatives,
we use variational derivatives.

\begin{myrule}[Factorial moment densities and Janossy densities]
The $k^{th}$-order density functions
can be determined from the p.g.fl. as follows.
\begin{align*}
m(x_1,\ldots,x_k)&=
\left.
{\delta^k 
G(h; \delta_{x_1},\ldots, \delta_{x_k}})\right|_{h=1}
\\
p(x_1,\ldots,x_k)&=
\left.
{\delta^k
G(h; \delta_{x_1},\ldots, \delta_{x_k}})
\right|_{h=0}
\end{align*}

\end{myrule}

\begin{example}[Poisson point process]
The p.g.fl. of a Poisson process is given by
\begin{align*}
G(h)=\exp(\lambda(f(h)-1)),
\end{align*}
where $f(h)=\int s(x)h(x)dx$.

\begin{align*}
\delta G(h; \eta)  =
\lambda f(\eta) \exp(\lambda(f(h)-1)).
\end{align*}
Setting $h=1$, we get
$
\lambda s(x)$.
More generally, we have
\begin{align*}
{\delta^r
G(h; \delta_{x_1},\ldots,\delta_{x_r}})=
\exp\left(\lambda(f(h)-1)\right)\prod_{i=1}^r \lambda s(x_i) 
\end{align*}
Then, setting $h=1$, we have
\begin{align*}
m(x_1,\ldots,x_r)&=
\prod_{i=1}^r \lambda s(x_i),
\end{align*}
or setting $h=0$, we get
\begin{align*}
p(x_1,\ldots,x_r)&=
\exp\left(-\lambda\right)
\prod_{i=1}^r \lambda s(x_i).
\end{align*}

\end{example}

\begin{myexercise}[Bernoulli point process]
Considering the p.g.fl. of the Bernoulli process
\begin{align*}
p_0+\int p_1(x)h(x)dx,
\end{align*}
find the first-order factorial moment density.
\end{myexercise}

\begin{myexercise}[Independent, identically distributed (i.i.d.) cluster point process]
Find the first-order factorial moment density
of the i.i.d. cluster process,
whose p.g.fl. is given by
\begin{align*}
G(h)=
p_0+
\sum_{k=1}^\infty
p_k\prod_{i=1}^k\int h(x_i)s(x_i)dx_i.
\end{align*}
\end{myexercise}

\newpage

\subsection{Joint probability generating functionals}

We now consider a motivating example
before extending to 
joint generating functionals.
\begin{example}
Consider
the joint functional
\begin{align*}
F(g,h)=
\int g(w)h(v) f(w,v) dw dv.
\end{align*}
First, let us take the functional 
derivative with respect to $g$
in the direction $\phi_g$,
\begin{align*}
\delta
F(g,h; \phi_g)=
\int \phi_g(w)h(v) f(w,v) dw dv.
\end{align*}
Now, take the functional derivative with
respect to $h$, i.e.
\begin{align*}
{\delta^2
}
 F(g,h;\phi_g;\phi_h)=
\int \phi_g(w)\phi_h(v) f(w,v) dw dv.
\end{align*}
Let us 
set $\phi_g=\delta_z$
and  $\phi_h=\delta_x$. 
Then 
\begin{align*}{\delta^2
}
F(g,h; \phi_g; \phi_h)=
\int \phi_g(w)\phi_h(v) f(w,v) dw dv
=f(z,x).
\end{align*}

Now let $f(z,x)=p(z,x)$
be a probability density function.
Note that,
according to the rules of conditional
probability,
$p(z,x)=p(z|x)p(x)$,
so that
\begin{align*}
\left.
{\delta^2
}
F(g,h;\phi_g;\phi_h)
\right|_{\phi_g=\delta_z,\phi_h=\delta_x}
=
p(z|x)p(x),
\end{align*}
and that
\begin{align*}
\left.
{\delta
}
F(g,1; \phi_g)
\right|_{\phi_g=\delta_z}
=
\int p(z|v)p(v)dv,
\end{align*}
which means that Bayes' rule can be written 
as
\begin{align*}
p(x|z)=
{
p(z|x)p(x)\over
\int p(z|v)p(v)dv
}
=
{
\left.
{\delta^2
}
F(g,h; \phi_g;\phi_h)
\right|_{\phi_g=\delta_z,\phi_h=\delta_x}
\over
\left.
{\delta
}
F(g,1;\phi_g) 
\right|_{\phi_g=\delta_z}
}.
\end{align*}

\end{example}

\subsubsection*{Bayesian estimation}

Now consider a joint functional \cite{Mahler_RPS_2007_3}
\begin{align*}
F(g,h)
=\sum_{n,m=0}^\infty
{1\over n!}
{1\over m!}
\left(\int \prod_{i=1}^n h(w_i) dw_i\right)
\left(\int \prod_{j=1}^m g(v_j) dv_j\right)
p(v_1,\ldots,v_m,w_1,\ldots,w_n).
\end{align*}
Then, clearly we have
\begin{align*}
\left.
{\delta^{n+m}
}F(g,h; \varphi_{z_1},\ldots, \varphi_{z_m};
\varphi_{x_1},\ldots,\varphi_{x_n})
\right|_{h=0,g=0}
=
p(z_1,\ldots,z_m,x_1,\ldots,x_n),
\end{align*}
when we set
$\varphi_{y}=\delta_{y}$.
Hence,
we can find 
\begin{align*}
p(x_1,\ldots,x_n|z_1,\ldots,z_m)
&=
{
p(z_1,\ldots,z_m,x_1,\ldots,x_n),
\over
p(z_1,\ldots,z_m),
}
\\
&=
{
\left.
{\delta^{n+m}
}F(g,h;
\varphi_{z_1}
,\ldots, \varphi_{z_m};
\varphi_{x_1},\ldots,\varphi_{x_n})
\right|_{h=0,g=0}
\over
\left.
{\delta^{m}
}F(g,1; \varphi_{z_1}
,\ldots, \varphi_{z_m})\right|_{g=0}
}
\end{align*}

It can often be convenient
to find the generating functional 
of $p(x_1,\ldots,x_n|z_1,\ldots,z_m)$,
which can be found with
\begin{align*}
G(h|z_1,\ldots,z_m)
&=
{
\left.
{\delta^{m}
}F(g,h;
\varphi_{z_1}
,\ldots, \varphi_{z_m})
\right|_{g=0}
\over
\left.
{\delta^{m}
}F(g,1;\varphi_{z_1}
,\ldots, \varphi_{z_m})\right|_{g=0}
}.
\end{align*}
We can determine the mean
from this 
with
\begin{align*}
D(x)=
\left.
{\delta}
G(h|z_1,\ldots,z_m;\varphi_x)\right|_{h=1}.
\end{align*}
This is known as the {\it Probability Hypothesis Density} or {\it PHD}.

\begin{myexercise}[The Probability
Hypothesis Density filter]

Let's look at the PHD filter update
(page 1173 of Mahler's 2003 AES paper),
and consider a joint functional of the form
\begin{align*}
F(g,h)=
G_\kappa(g)G_p(hG_L(g|\cdot)),
\end{align*}
where 
\begin{align*}
G_\kappa(g)=\exp(\lambda(\int c(z)g(z)dz -1)),
\end{align*}
 is the Poisson p.g.fl. 
for false alarms
where $\lambda$
is the average number of false alarms;
\begin{align*}
G_p(h)=\exp(\mu(\int s(x)h(x)dx-1))
\end{align*}
is the prior p.g.fl.,
where $\mu$ is the average number of targets,
each of which distributed according to $s(x)$;
and 
\begin{align*}
G_L(g|x)=1-p_D(x)+p_D(x)\int g(z)f(z|x)dz
\end{align*}
is the Bernoulli detection
process for each target
where the probability of detecting each target
is $p_D(x)$
and the 
likelihood of observing
measurement $x$
conditioned on $x$ is given by
$f(z|x)$.
In this case, we have the joint p.g.fl.
\begin{align*}
F(g,h)=
\exp
\left(
\lambda(\int c(z)g(z)dz-1)
+\mu(\int s(x)h(x) G_L(g|x)dx -1)
\right).
\end{align*}

To find the PHD filter update, 
firstly find the functional derivative of $F(g,h)$
with respect to
$g$ in the direction $\varphi_{z_1}=\delta_{z_1}$.
Now differentiate again in the direction
$\varphi_{z_2}=\delta_{z_2}$.
Generalise to find the $m^{th}$ derivative of $F(g,h)$
in directions $\delta_{z_1}, \delta_{z_2},\ldots, \delta_{z_m}$.
Now determine the 
updated p.g.fl. $G(h|z_1,\ldots,z_m)$.
Differentiate again to obtain 
the updated PHD, $D(x|z_1,\ldots,z_m)$
and show that it is equal to
\begin{align*}
D(x|z_1,\ldots,z_m)=
\mu s(x)(1-p_D(x))
+\sum_{z\in Z}
{p_D(x)f(z|x)\mu s(x)
\over
\lambda c(z)+\mu \int p_D(x) s(x) f(z|x) dx
}.
\end{align*}
\end{myexercise}

\begin{myexercise}
Suppose that $G_p(h)=1-p+p f(h)$, where $f(h)=\int s(x)h(x)dx$,
and find the updated PHD \cite{Mahler_RPS_2003}.
\end{myexercise}

\end{document}